%% file: BDT.tex
\documentclass[11pt]{article}
\usepackage[pagebackref,colorlinks=true,urlcolor=red,linkcolor=red,citecolor=red]{hyperref}

\usepackage{epsfig,graphics, graphicx}
\usepackage{subcaption, comment, bm}
\usepackage[percent]{overpic}

\usepackage{amssymb,bm}
\usepackage{amsthm}
\usepackage{color}
\usepackage[]{amsmath}
\usepackage[]{amsfonts}
\usepackage[]{fancyhdr}
\usepackage[]{graphicx,wrapfig}
\usepackage{enumitem}
\usepackage{array}

\graphicspath{{/EPSF/}{../figures/}{figures/}}

\input preamble

\usepackage{graphicx}
\usepackage[dvipsnames]{xcolor}

\newcommand{\Ic}{\mathcal{I}}
\newcommand{\Jc}{\mathcal{J}}

\newcommand{\Lc}{\mathcal{L}}

\newcommand{\Sc}{\mathcal{S}}
\newcommand{\Tc}{\mathcal{T}}

\newcommand{\Xc}{\mathcal{X}}

\newcommand{\Rb}{\mathbb{R}}

\newcommand{\vP}{\textbf{\textit{P}}}
\newcommand{\vf}{\textbf{\textit{f}}}
\newcommand{\vx}{\textbf{\textit{x}}}
\newcommand{\vy}{\textbf{\textit{y}}}
\newcommand{\vz}{\textbf{\textit{z}}}
\newcommand{\vw}{\textbf{\textit{w}}}
\newcommand{\vu}{\textbf{\textit{u}}}
\newcommand{\vv}{\textbf{\textit{v}}}
\newcommand{\vn}{\textbf{\textit{n}}}
\newcommand{\vgamma}{\pmb{\gamma}}
\newcommand{\vpsi}{\pmb{\psi}}
\newcommand{\vzeta}{\pmb{\zeta}}
\renewcommand{\O}{\Omega}

\usepackage{amsthm}
\usepackage{amsfonts}
\newcommand{\Beq}{\begin{equation}}
\newcommand{\Eeq}{\end{equation}}
\newcommand{\beq}{\begin{equation*}}
\newcommand{\eeq}{\end{equation*}}
\newcommand{\bal}{\begin{align}}
\newcommand{\eal}{\end{align}}

\newtheorem*{theor}{Theorem}
\newtheorem{thr}{Theorem}
\newtheorem{defn}{Definition}
\newtheorem{rem}{Remark}
\newtheorem{cor}{Corollary}

\title{\vspace{-1cm} Generalized V-line transforms in 2D vector tomography}

\author{Gaik Ambartsoumian\thanks{Department of Mathematics, University of Texas at Arlington, Arlington, TX, United States of America.  gambarts@uta.edu}\and Mohammad Javad Latifi Jebelli\thanks{Department of Mathematics, University of Arizona, Tucson, AZ, United States of America. mjlatifi@math.arizona.edu} \and Rohit Kumar Mishra\thanks{Department of Mathematics, University of Texas at Arlington, Arlington, TX, United States of America. rohit.mishra@uta.edu}}

\begin{document}
\date{}
\maketitle
\begin{abstract}
We study the inverse problem of recovering a vector field in $\Rb^2$ from a set of new generalized $V$-line transforms in three different ways. First, we introduce the longitudinal and transverse $V$-line transforms for vector fields in $\Rb^2$. We then give an explicit characterisation of their respective kernels and show that they are complements of each other. We prove invertibility of each transform modulo their kernels and combine them to reconstruct explicitly the full vector field. In the second method, we combine the longitudinal and transverse V-line transforms with their corresponding first moment transforms and recover the full vector field from either pair. We show that the available data in each of these setups can be used to derive the signed V-line transform of both scalar component of the vector field, and use the known inversion of the latter. The final major result of this paper is the derivation of an exact closed form formula for reconstruction of the full vector field in $\mathbb{R}^2$ from its star transform with weights. We solve this problem by relating the star transform of the vector field to the ordinary Radon transform of the scalar components of the field.
\end{abstract}
\vspace{-5mm}
\section{Introduction}\label{Introduction}

The primary task of integral geometry is reconstructing a scalar function or a vector field (or more generally a tensor field) from some kind of integral transform data. These types of reconstruction problems are often crucial parts of various non-invasive imaging techniques with applications in medicine, seismology, oceanography and  many other areas.

A typical integral geometry problem can be formulated as follows. What information about a tensor field of rank $m$ can be recovered from its longitudinal ray transform (also known as ray transform)? It has been shown by several authors that a scalar function (corresponding to the case $m$=0) can be reconstructed uniquely from the knowledge of its ray transform, e.g. see \cite{Nattrer_Wuebbeling}. For $m \geq 1$, this transform has a non-trivial kernel, which makes the full recovery of a tensor field impossible, when using only the ray transform data. In this case, only the solenoidal part of a tensor field can recovered. The latter problem has been studied  in various settings by multiple authors, e.g. see \cite{Denisjuk1994,Denisjuk_Paper,Katsevich2006,Katsevich2013,Venky_and_Rohit,Monard2016a,Palamodov2009, PSU2012, Schuster2000,Sharafutdinov_Book,Sharafutdinov2007,Tuy1983} and the references therein.

The non-injectivity of the ray transform raises a natural question: what kind of additional data is needed for full reconstruction of tensor fields? In this context, an injectivity result has been presented utilizing the so-called integral moment transforms over symmetric $m$-tensor fields in $\Rb^n$, see \cite{Sharafutdinov_Generalized_Tensor_Fields}. We also point out some recent related works for the invertibilty of these generalized transforms \cite{Anuj_Rohit, Wongsason_2020,Venky_Suman_Manna,Venky_Suman_Manna2,Francois_Rohit_Venky, R_K_Mishra,Rohit_Suman_2020}.

Another approach to reconstruct the full tensor field is to work with the transverse ray transform (TRT) instead of (or in addition to) the longitudinal ray transform (LRT).  For $n =2$, TRT and LRT provide equivalent information, up to a linear transformation of the tensor field. In particular, TRT also has a non-trivial kernel, making the full recovery just from TRT impossible. However, one can combine the data from both transforms (TRT and LRT) in 2D to recover a vector field completely \cite{Derevtsov2}. {For the recovery of tensor fields, one needs to work with mixed ray transforms, which are natural generalizations of TRT and LRT, e.g. see \cite{Derevtsov3,Hoop2019}}. In contrast to the 2D case, when $n \geq 3$ it is known that a symmetric $m$-tensor field is completely determined by its transverse ray transform, see \cite{Anuj_TRT,Derevtsov3, Holman2013, Novikov_Sharafutdinov, Sharafutdinov_TRT_2008}. In addition to these injectivity results, there are also various reconstruction results for TRT in different settings ($n \geq 3$), see  \cite{Rohit_Griemaier,VRS,Wongsason_2018} and reference therein.

%--------------------------------------------
\vspace{2mm}
In this article, we consider a full reconstruction of a vector field in $\mathbb{R}^2$ using a new set of integral transforms (see Tables \ref{Table1} and \ref{Table2}). These operators are analogous to the ray transforms discussed above, but instead of integrating along straight lines they use V-shaped paths of integration, called V-lines.

%in three different ways: simultaneously using its longitudinal V-line transform ($\mathcal{L}$) and transverse V-line transform ($\mathcal{T}$) (see Definitions \ref{def:definition of V line Doppler transform} and \ref{def:definition of transverse V line Doppler transform}), using its $\mathcal{L}$ together with the first moment longitudinal V-line transform ($\mathcal{I}$) (see Definition \ref{def:definition of first V line moment  transform}) or alternatively using its $\mathcal{T}$ together with first moment transverse V-line transform ($\mathcal{J}$) (see Definition \ref{def:definition of first V line transverse moment  transform}), and using its star transform (see Definition \ref{def:star}).

The V-line transform (often also called broken ray transform) for \textit{scalar functions} in $\mathbb{R}^2$ maps a function to its integrals along piecewise linear trajectories, which consist of two rays emanating from a common vertex. Two distinct classes of V-line transforms with some generalizations have been studied by various authors in recent past. The first class includes V-lines (and cones in higher dimensions) that have a vertex on the boundary of the image domain, i.e. outside of the support of the image function (see the review article \cite{Terzioglu_2018} and the references there). These transforms often appear in image reconstruction problems using Compton cameras. The second class includes V-lines (as well as stars and cones) that have a vertex inside the image domain, and they appear in relation to  single scattering tomography \cite{Gaik_2012, Amb2, Ambartsoumian_2019, Amb_Lat_star, Gaik_Moon, Amb-Roy, Florescu_Markel_Schotland_2011, Florescu_Markel_Schotland_2009, Florescu3, Florescu4, Rim_Gaik_2014, Katsevich_2013, Kats-Kryl-15, Sherson-15, Terzioglu_2015, Walker_2019, ZSM-star-14}. The V-line transforms of \textit{vector fields} discussed in this paper are a natural generalization of the second class of V-lines discussed above.

%------------------------------------------------

\begin{table}
\centering
 \begin{tabular}{||c l c ||}
 \hline
 Symbol & Name & Definition  \\ [0.5ex]
 \hline\hline
 $\mathcal{L}$ & Longitudinal V-line transform & \ref{def:definition of V line Doppler transform}  \\
 \hline
  $\mathcal{T}$ & Transverse V-line transform & \ref{def:definition of transverse V line Doppler transform} \\
 \hline
 $\mathcal{I}$ & First moment longitudinal V-line transform & \ref{def:definition of first V line moment  transform}  \\

 \hline
 $\mathcal{J}$ & First moment transverse V-line transform &  \ref{def:definition of first V line transverse moment  transform}  \\
 \hline
 $\mathcal{S}$ & Vector-valued star transform & \ref{def:star_d}\\  %[1ex]
 \hline
\end{tabular}
 \caption{A list of integral operators discussed in the paper.}
\label{Table1}
\end{table}
%\hspace{2cm}

\begin{table}
\centering
 \begin{tabular}{||l c ||}
 \hline
Reconstruction of $\vf$ from: & Theorem   \\ [0.5ex]
 \hline\hline
   knowledge of $\mathcal{L}\vf$ and $\mathcal{T}\vf$ & \ref{th: Inversion of B}, \ref{th: Inversion of T}   \\
 \hline
  knowledge of $\mathcal{L}\vf$ and $\mathcal{I}\vf$ & \ref{th:inversion using moment transform}   \\
 \hline
  knowledge of $\mathcal{T}\vf$ and $\mathcal{J}\vf$ & \ref{th:inversion using transverse moment transform}  \\
 \hline
 knowledge of $\mathcal{S}\vf$ &  \ref{th:star}  \\ %[1ex]
 \hline
\end{tabular}
\caption{A list of reconstructions provided in the paper.}
\label{Table2}
\end{table}

The rest of the paper is organized as follows. In Section \ref{sec:def} we introduce the notations and define the operators used in this article. In Section \ref{sec:LVT-TVT} we state the main results about the kernels of $\mathcal{L}$ and $\mathcal{T}$, as well as the relations of these transforms correspondingly to the curl and divergence of the vector field enabling its full recovery. In Section \ref{sec:Proofs} we provide the proofs of theorems stated in the previous section. Section \ref{sec:LVT-moment} describes the method of recovering the full vector field from either one of the pairs: $\mathcal{L}$ with its first moment $\mathcal{I}$, and $\mathcal{T}$ with its first moment $\mathcal{J}$. Section \ref{sec:star} presents an exact closed form formula for recovering the full vector field from its star transform. We finish the paper with some additional remarks listed in Section \ref{sec:remarks} and acknowledgements in Section \ref{sec:acknowledge}.
%%%%%%%%%%%%%%%%%%%%%%%%%%%%%%%%%%%%%%%%%%%%%%%%%%%%%
\section{Definitions and notations}\label{sec:def}

In this section we introduce the notations and define the operators used in the article. Throughout the paper, we use bold font letters to denote vectors in $\mathbb{R}^2$ (e.g. $\vx$, $\textbf{\textit{u}}$, $\textbf{\textit{v}}$, $\textbf{\textit{f}}$, etc), and regular font letters to denote scalars (e.g. $t$, $h$, $f_i$, etc). We denote by $\vx\cdot\vy$ the usual dot product between vectors $\vx$ and $\vy$.

For a scalar function $V(x_1, x_2)$  and a vector field $\vf =(f_1,f_2)$, we use the notations
\begin{align}\label{eq: definition of div and curl}
\nabla V := \left(\frac{\partial V}{\partial x_1}, \frac{\partial V}{\partial x_2}\right), \quad   \operatorname{div}\vf :=  \frac{\partial f_1}{\partial x_1}+ \frac{\partial f_2}{\partial x_2},\quad \mbox{ and } \quad  \operatorname{curl} \vf :=  \frac{\partial f_2}{\partial x_1}- \frac{\partial f_1}{\partial x_2}.
\end{align}
{The operators $\nabla$ and $\operatorname{div}$ are the classical gradient and divergence operators respectively. The operator $\operatorname{curl}$ defined above is essentially the exterior derivative on 2-dimensional manifolds (e.g. see \cite[Chapter 14]{Lee_Book}). However, it is customary to call that operator $\operatorname{curl}$ in 2D
(e.g. see \cite{Derevtsov2, Derevtsov3,Sparr1995}).}
%The goal of this paper is to recover a vector field from the knowledge of two integral transforms, namely, the V-line Doppler transform and the V-line transverse ray transform. Each of these transforms has a non-trivial kernel, which makes the recovery of the full vector field from only one of them impossible, motivating the simultaneous use of data from both transforms.

%\begin{figure}{}
%\centering
%\includegraphics[height=3.8cm]{v-line-def.png}\caption{V-line with vertex at $\vx$ and ray directions $\vu$ and $\vv$}
%\label{fig:vline-def}
%\end{figure}

\begin{figure}[ht]
\begin{center}
\begin{subfigure}{.45\textwidth}
\centering
 \includegraphics[height=3.8cm]{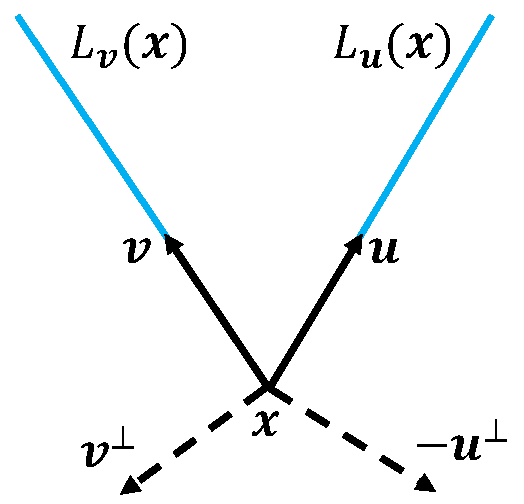}
 \caption{A V-line with vertex at $\vx$, ray directions $\vu$, $\vv$ and outward normals $-\vu^\perp$, $\vv^\perp$.}
  \label{fig1a}
\end{subfigure} \qquad
\begin{subfigure}{.45\textwidth}
\centering
\includegraphics[height=3.8cm]{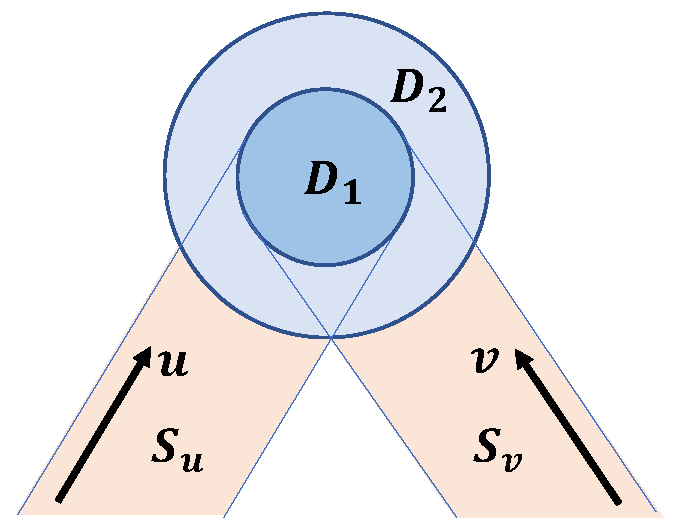}
\caption{A sketch of the compact support of  $\vf$ and the unbounded support of $\Lc\vf$, $\Tc\vf$, $\Ic\vf$, $\Jc\vf$.} \label{fig1b}
\end{subfigure}
\end{center}
\vspace{-5mm}
\caption{}
\label{fig1}
\end{figure}

Let  $\vu$ and $\vv$ be two linearly independent unit vectors in $\mathbb{R}^2$. For $\vx \in \mathbb{R}^2$, the rays emanating from $\vx$ in directions $\vu$ and $\vv$ are denoted by  $L_{\vu}(\vx)$ and $L_{\vv}(\vx)$ respectively, i.e.
$$ L_{\vu}(\textbf{\textit{x}}) = \left\{\textbf{\textit{x}} +t \textbf{\textit{u}}: 0 \leq t < \infty\right\} \quad \mbox{ and } \quad  L_{\vv}(\textbf{\textit{x}}) = \left\{\textbf{\textit{x}} +t \textbf{\textit{v}}: 0 \leq t < \infty\right\}.$$
A V-line with vertex $\textbf{\textit{x}}$ is the union of rays $L_{\vu}(\textbf{\textit{x}})$ and $L_{\vv}(\textbf{\textit{x}})$. For the rest of the article we will assume that  $\textbf{\textit{u}}$ and $\textbf{\textit{v}}$ are fixed, i.e. all V-lines have the same ray directions and can be parametrized simply by the coordinates $\textbf{\textit{x}}$ of the vertex (see Figure \ref{fig1a}).

\begin{defn}
 The \textbf{divergent beam transform} $\mathcal{X}_{\vu}$ of function $h$ at $\vx\in \mathbb{R}^2$ in the direction $\vu$ is defined as:
 \begin{equation}\label{def:DivBeam}
   \mathcal{X}_{\vu}h(\vx) =  \int_{0}^{\infty} h(\vx+t \vu)\,dt.
 \end{equation}
 \end{defn}
 \noindent The directional derivative of a function in the direction $\vu$ is denoted by $D_{\vu}$, i.e.
\begin{equation}
D_{\vu} h = \vu \cdot \nabla h.
%u_1 \frac{\partial h}{\partial x_1} + u_2 \frac{\partial h}{\partial x_2}
\end{equation}
One can similarly define the divergent beam transform $\mathcal{X}_{\vv}$ and directional derivative $D_{\vv}$.\vspace{2mm}\\
{For  a compactly supported  function $h$, we observe the inverse relations (following from the fundamental theorem of calculus) between the operators $\Xc_{\vu}$ and $D_{\vu}$ (and similarly between $\Xc_{\vv}$ and $D_{\vv}$)
\begin{align}\label{eq:intertwining relations}
  \Xc_{\vu}\left(D_{\vu} h\right) (\vx) = -h(\vx) \quad \mbox{ and  } \quad   D_{\vu}\left(\Xc_{\vu} h\right) (\vx) = -h(\vx).
\end{align}}
%%%%%%%%%%%%%%%%%%%%%%%%%%%%%%%%%%%%%%%%%%%%%%%%
\begin{comment}
Before defining the transforms of interest, we would like to introduce notations for some integral and differential operators. For any given function $h$,
\begin{itemize}
\item $\mathcal{X}_{\vu}(h)$ and $\mathcal{X}_{\vv}(h)$ denote  the integrals of $h$ along $L_{\vu}(\vx)$ and $L_{\vv}(\vx)$  respectively. More precisely,
$$ \mathcal{X}_{\vu}h(\vx) = \int_{L_{\vu}(\vx)}h\ d l, \quad \mbox{ and } \quad \mathcal{X}_{\vv}h(\vx) = \int_{L_{\vv}(\vx)} h\ d l,$$
where $d l$ is the standard Lebesgue measure on the line.
\item The directional derivatives in the directions $\vu$ and $\vv$ will be denoted by $D_{\vu}$ and $D_{\vv}$ respectively, i.e.
$$ D_{\vu} h = u_1 \frac{\partial h}{\partial x_1} + u_2 \frac{\partial h}{\partial x_2} \quad \mbox{ and } \quad D_{\vv} h = v_1 \frac{\partial h}{\partial x_1} + v_2 \frac{\partial h}{\partial x_2}.$$
\end{itemize}
\end{comment}
%%%%%%%%%%%%%%%%%%%%%%%%%%%%%%%%%%%%%%%%%%%%%

%Now, we are ready to define the integral transforms studied in this manuscript.

The goal of this paper is to recover a vector field from the knowledge of its various integral transforms, namely: $\mathcal{L}$, $\mathcal{T}$, $\mathcal{I}$, $\mathcal{J}$ and the star transform $\mathcal{S}$. These transforms are defined in analogy with the corresponding ray transforms of vector fields in $\mathbb{R}^2$, substituting the straight line trajectory of integration of the latter with a $V$-line or star trajectory for the former.
%Each of these transforms has a non-trivial kernel, which makes the recovery of the full vector field from only one of them impossible, motivating the simultaneous use of data from both transforms.
In applications, V-lines correspond to flight paths of particles that scatter at some point in the medium.
%Physically, the direction of traveling along a V-line corresponds to the flight path of a scattering particle.
One can imagine a particle starting from a point at infinity traveling along direction $-\vu$ to $\vx$, where scattering happens, after which the particle goes to infinity in the direction of $\vv$. This discussion motivates the following:
\begin{defn} \label{def:definition of V line Doppler transform} Let $\vf = (f_1, f_2)$ be a vector field in $\mathbb{R}^2$ with components $f_i \in C^2_c(\Rb^2)$ for $i =1, 2$. The \textbf{longitudinal V-line transform} of $\textbf{\textit{f}}$ is defined as
\begin{align}\label{eq:def V-line transform}
\mathcal{L}_{\vu, \vv}\, \vf\  = -\mathcal{X}_{\vu}  \left(\vf \cdot \vu\right)  + \mathcal{X}_{\vv}  \left(\vf \cdot \vv\right).
% \mathcal{L}_{\vu, \vv}\, \vf\ (\vx) =  -\int_0^\infty \left\langle \vf\, (\vx +t \vu), \vu \right\rangle dt + \int_0^\infty \left\langle \vf\, (\vx +t \vv), \vv \right\rangle dt,
\end{align}
% i.e.
% $$
% \mathcal{L}_{\vu, \vv}\, \vf\  = -\mathcal{X}_{\vu}  \left(\vf \cdot \vu\right)  + \mathcal{X}_{\vv}  \left(\vf \cdot \vv\right),
% $$
%where $\vf \cdot \vu$ denotes the usual dot product of $\vf$ and $\vu$.
\end{defn}
%\begin{wrapfigure}{r}{5cm}
%\includegraphics[width=5cm]{NormaL_{\vv}ectors}\caption{Normal vectors}
%\end{wrapfigure}
%\begin{remark}
%Physically the direction of traveling along the V-line corresponds to the flight path of a scattering particle. This is why the inner product of the unknown vector field $\vf$ in the first integral is taken with $-\vu$ %instead of $\vu$.
%\end{remark}

To define the second integral transform of interest, we need to make a choice for the normal unit vector corresponding to each branch of the V-line. We define the vector $\perp$ operation by $(x_1,x_2)^\perp = (-x_2,x_1) $.

%\begin{figure}
%\begin{center}
% \includegraphics[width=6cm]{Normalvectors.eps}\caption{Choice of normal vectors}
%\end{center}
%\end{figure}

%--------------------------------------------------

\begin{defn}\label{def:definition of transverse V line Doppler transform}
Let $\vf = (f_1, f_2)$ be a vector field in $\mathbb{R}^2$ with components $f_i \in C^2_c(\Rb^2)$ for $i =1, 2$. The \textbf{transverse V-line transform} of $\textbf{\textit{f}}$ is defined as
\begin{align}\label{eq:def transverse V-line transform}
\mathcal{T}_{\vu, \vv}\, \vf\  = - \mathcal{X}_{\vu}  \left(\vf \cdot \vu^\perp\right) + \mathcal{X}_{\vv}  \left(\vf \cdot \vv^\perp\right).
\end{align}
%which clearly satisfies $\mathcal{T}_{\vu, \vv}\, \vf = -\mathcal{L}_{\vu, \vv}\, \vf^\perp$.
\end{defn}

%--------------------------------------------------

The orientation of normal vectors is chosen towards the same side of the path of the scattering particle.
Hence, in the definition above the inner product of the unknown vector field is taken with the outward unit normal of the V-line at each point (see Figure \ref{fig1a}).

%--------------------------------------------------
%\begin{rem}
% It is easy to verify that $\mathcal{T}_{\vu,\vv} \vf = -\mathcal{L}_{\vu,\vv} \vf^{\,\perp}$.\\
% \end{rem}

%-----------------------------------------------

\begin{defn}
The \textit{\textbf{first moment  divergent beam transform}} of a function $h$ in the direction $\vu$ is defined as follows  $$\Xc^1_{\vu} h(\vx) =  \int_0^\infty h(\vx + t \vu)\, t\,  dt.  $$
\end{defn}
\noindent Similarly
$\Xc^1_{\vv} h(\vx) =  \int_0^\infty  h(\vx + t \vv)\,t\, dt.$
\begin{defn} \label{def:definition of first V line moment  transform} Let $\vf = (f_1, f_2)$ be a vector field in $\mathbb{R}^2$ with components $f_i \in C^2_c(\Rb^2)$ for $i =1, 2$. The \textbf{first moment longitudinal V-line transform} of $\textbf{\textit{f}}$ is defined as
\begin{align}\label{eq:first V line moment  transform}
\mathcal{I}_{\vu, \vv}\, \vf\ (\vx) &=  - \Xc_{\vu}^1 (\vf \cdot \vu) + \Xc_{\vv}^1(\vf\cdot \vv).
\end{align}
\end{defn}

\begin{defn} \label{def:definition of first V line transverse moment  transform} Let $\vf = (f_1, f_2)$ be a vector field in $\mathbb{R}^2$ with components $f_i \in C^2_c(\Rb^2)$ for $i =1, 2$. The \textbf{first moment transverse V-line transform} of $\textbf{\textit{f}}$ is defined as
\begin{align}\label{eq:first V line moment  transform2}
\mathcal{J}_{\vu, \vv}\, \vf\ (\vx) &=  - \Xc_{\vu}^1 \left(\vf \cdot \vu^\perp \right) + \Xc_{\vv}^1\left(\vf\cdot \vv^\perp \right).
\end{align}
\end{defn}
\vspace{2mm}

%-----------------------------------------------
\begin{rem}
It is easy to verify that $\mathcal{T}_{\vu,\vv} \vf = -\mathcal{L}_{\vu,\vv} \vf^{\,\perp}$ and $\Jc_{\vu, \vv} \vf = - \Ic_{\vu, \vv} \vf^{\,\perp}$.
\end{rem}

%-----------------------------------------------
\vspace{2mm}
\begin{rem}
Throughout the paper we assume that the linearly independent unit vectors $\vu$ and $\vv$ are fixed. Hence, to simplify the notations we will drop the indices $\vu,\vv$ and refer to $\mathcal{T}_{\textbf{\textit{u}}, \textbf{\textit{v}}}$, $\mathcal{L}_{\textbf{\textit{u}},  \textbf{\textit{v}}}$, $\mathcal{I}_{\vu, \vv}$ and $\Jc_{\vu, \vv}$ simply as $\mathcal{T}$, $\mathcal{L}$, $\mathcal{I}$ and $\Jc$.
\end{rem}

%-------------------------------------------------

\noindent {\textbf{The support of $\vf$ and related restrictions of its transforms}}

\noindent Let us assume that  $\operatorname{supp}\vf\subseteq D_1$, where $D_1$ is an open disc of radius $r_1$ centered at the origin. {Let $D_2$ be the smallest disc (of some finite radius $r_2>r_1$) centered at the origin, such that only one ray of any V-line with a vertex outside of $D_2$ intersects $D_1$ (see Figure \ref{fig1b}).} Then $\Lc\vf$, $\Tc\vf$, $\Ic\vf$ and $\Jc \vf$ are supported inside an unbounded domain $D_2\cup S_{\vu} \cup S_{\vv}$, where
%$D_2$ is a disc of some finite radius $r_2>r_1$ centered at the origin, while
$S_{\vu}$ and $S_{\vv}$ are semi-infinite strips (outside of $D_2$) in the direction of $\vu$ and $\vv$ correspondingly (see Figure \ref{fig1b}). It is easy to notice that all three transforms $\Lc\vf$, $\Tc\vf$, $\Ic\vf$ and $\Jc \vf$ are constant along the directions of rays $\vu$ and $\vv$ inside the corresponding strips $S_{\vu}$ and $S_{\vv}$. In other words, the restrictions of $\Lc\vf$, $\Tc\vf$,  $\Ic\vf$ and $\Jc \vf$ to $\overline{D_2}$ completely define them in $\mathbb{R}^2$.

\begin{rem}
Throughout the paper we assume that { the vector field $\vf$ is supported in $D_1$ and the transforms} $\Lc\vf\,(\vx)$, $\Tc\vf\,(\vx)$, $\Ic\vf\,{(\vx)}$ and $\Jc\vf\,{(\vx)}$ are known for all $\vx\in \overline{D_2}$.
\end{rem}

%%%%%%%%%%%%%%%%%%%%%%%%%%%%%%%%%%%%%%%%%%%%%%%%%%%%%%%%
\section{Full field recovery using longitudinal and transverse VLT}\label{sec:LVT-TVT}
%\subsection{Statement of main results}
The following two relations can be obtained by a simple calculation:
\begin{align}
\Delta f_1 &= \frac{\partial}{\partial x_1}\operatorname{div}  \vf -\frac{\partial}{\partial x_2}\operatorname{curl}  \vf \label{eq:Laplace of components of f1},\\
\Delta f_2 &= \frac{\partial}{\partial x_2}\operatorname{div}  \vf +\frac{\partial}{\partial x_1}\operatorname{curl}  \vf \label{eq:Laplace of components of f2}.
\end{align}
Therefore the Laplacian of each component of a vector field $\vf$ can be computed explicitly if one knows $\operatorname{div}\vf$ and $\operatorname{curl}\vf$. {The knowledge of  Laplacians of each component of $\vf$ allows an explicit reconstruction of $\vf$ with the help of Green's function on the disc $D_1$}. Hence, we have the following

\begin{rem}
To recover a compactly supported vector field $\vf$ {explicitly}, one only needs to reconstruct $\operatorname{div}\vf$ and $\operatorname{curl}\vf$ from the integral transforms under consideration.
\end{rem}

The following two theorems show that there is a non-trivial kernel for each of the integral transforms $\Lc$ and $\Tc$. Moreover, the theorems explicitly characterize those kernels.
\begin{thr}\label{th: kernel of B}
The kernel of longitudinal V-line transform $\Lc$ is the set of all  potential vector fields $\vf$. In other words, {if $\vf$ is a vector field in $\Rb^2$ with components in $C_c^2(D_1)$, then}
\begin{align*}
\Lc \vf \equiv0 \quad \mbox{if and only if\ \ } \vf = \nabla V,  \ \mbox{ for some scalar function } V.
\end{align*}
\end{thr}
\noindent
One can easily check that all potential vector fields $\vf = \nabla V$ are curl-free (i.e. curl $\vf =0$) and vice versa. Thus, from the above theorem we conclude that all curl-free vector fields are in the kernel of $\Lc$.
\begin{thr}\label{th: kernel of T}
The kernel of transverse V-line transform $\Tc$ is the set of all divergence-free vector fields $\vf$. In other words, {if $\vf$ is a vector field in $\Rb^2$ with components in $C_c^2(D_1)$, then}
\begin{align*}
\Tc \vf \equiv 0 \quad \mbox{if and only if }\ \  \operatorname{div} \vf =0.
\end{align*}
\end{thr}
Before moving on, we would like to recount here a crucial and well known theorem, which states that any vector field (with some boundary condition) can be decomposed uniquely into a divergence-free part and a curl-free part. The following decomposition result is true in more general settings, e.g. in arbitrary dimensions, as well as for tensor fields. But for our needs, it is sufficient to consider the statement just for vector fields in $\mathbb{R}^2$.
\begin{theor}[Theorem 3.3.2, \cite{Sharafutdinov_Book}]\label{th: decomposition of tensor field on manifolds}
Let $\O$ be a bounded domain in $\Rb^2$ and $\vf$ be a vector field,  whose support is contained in $\O$. Then there exist a uniquely determined vector field $\textup{\vf}^s$ and a uniquely determined scalar function $V$ satisfying
	\begin{align}\label{eq:decomposition}
	\vf = \textup{\vf}^s + \nabla V \quad  \text{ with } \quad \operatorname{div} \textup{\vf}^s =0 \quad \text{ and } \quad V|_{\partial \O}=0.
	\end{align}

%	Let $M$ be a compact Riemannian manifold with boundary. Let $k\geq1$ be an integer. Then for any vector field $\vf$ with components in $H^k(M)$ there exist uniquely determined vector field $\vf^{\ s}$ with components in $H^k(M)$ and a scalar function $V\in H^{k+1}(M)$ such that
%	\begin{align}
%	\vf = \textup{\vf}^s + \nabla V \quad  \text{ with } \quad \textup{div } \textup{\vf}^s =0 \quad \text{ and } \quad V|_{\partial M}=0.
%	\end{align}
\end{theor}
The fields $\vf^s$ and   $\nabla V$ are known as the solenoidal part (divergence-free part) and the potential part (curl-free part) of $\vf$ respectively.\\
	
From Theorems \ref{th: kernel of B} and \ref{th: kernel of T} we see that the solenoidal part $\vf^s$ and the potential part $\nabla V$ of $\vf$ are always in the kernel of $\Tc$ and $\Lc$ respectively. Hence it is impossible to reconstruct the full vector field just from the knowledge of only one transform $(\Lc$ or $\Tc$). Also, observe from the above decomposition that
\begin{align*}
\operatorname{curl}  \vf = \operatorname{curl}  \vf^s \quad \mbox{ and } \quad \operatorname{div}  \vf =  \Delta V.
\end{align*}
This implies that the problem of recovering $\operatorname{div}\vf$ and $\operatorname{curl}\vf$  is  reduced to the determination of  $\Delta V$ and $\operatorname{curl}\vf^s$. Our next two theorems state that it is indeed possible to reconstruct $\Delta V$ and $\operatorname{curl}\vf^s$ explicitly from the knowledge of $\Tc \vf$ and $\Lc \vf$ respectively.
\begin{thr}\label{th: Inversion of B}
Let $\vf$ be a vector field in $\Rb^2$ with components in $C_c^2({D_1})$. Then
$\operatorname{curl} \vf$ can be recovered from $\Lc \vf$ as follows:
\begin{align}\label{eq:curl-explicit}
\operatorname{curl} \vf &= \frac{1}{\det(\vv, \vu)}D_{\vu}D_{\vv}\, \Lc \vf.
\end{align}
In particular, this implies that operator $\Lc$ is invertible over compactly supported divergence-free vector fields.
\end{thr}
\begin{thr}\label{th: Inversion of T}
Let $\vf$ be a vector field   in $\Rb^2$ with components in $C_c^2({D_1})$. Then $\operatorname{div} \vf$ can be recovered from  $\Tc \vf$ as follows:
\begin{align}\label{eq:div-explicit}\operatorname{div} \vf &= -\frac{1}{\det(\vv, \vu)}D_{\vu}D_{\vv} \Tc \vf. \end{align}
 In particular, this implies that operator $\Tc$ is invertible over compactly supported curl-free vector fields.
\end{thr}
\begin{rem}
The quantity appearing in the denominator of expressions for $\operatorname{curl} \vf$ and $\operatorname{div} \vf$ is not zero, since $\vu$ and $\vv$ are linearly independent. In other words,
$$
\det(\vv, \vu)=v_1u_2-u_1v_2=\vu\cdot\vv^{\perp}\ne0.
$$
\end{rem}
% \rohit{I feel we should remove both these remarks as it is already included in the statement of theorem ``in particular parts of statements". Basically, we have that the }

In some cases one may be interested in an unknown scalar potential $V$ supported in $D_1$, while only having measurements $\Tc\vf$ of its gradient $\vf=\nabla V$. Since $\operatorname{div} \vf = \Delta V$,  as a consequence of Theorem \ref{th: Inversion of T} we can recover the scalar function  $V$ explicitly by solving the following Dirichlet problem for the Poisson equation:
\begin{align*}
    \left\{\begin{array}{lll}
    &\Delta V(\vx)  =    \displaystyle -\frac{1}{\det(\vv, \vu) }D_{\vu}D_{\vv} \Tc \vf\,(\vx)&  \mbox{\textup{ in }} D_1,  \\
      &V(\vx) = 0  &\mbox{ \textup{on} } \partial D_1.
    \end{array}\right.
\end{align*}

Similarly, one may be interested in a compactly supported scalar function $W$, when the measurements $\Lc\vf$ are available only for
%its (distorted) gradient
$ \vf = (\nabla W)^\perp = \left( -\frac{\partial W}{\partial x_2}, \frac{\partial W}{\partial x_1}\right)$. In such cases,  one may use the relation $\operatorname{curl} \vf = \Delta W$ to get $W$ by solving the following Dirichlet boundary value problem  \begin{align*}
    \left\{\begin{array}{lll}
    &\Delta W(\vx)  =    \displaystyle \frac{1}{\det(\vv, \vu) }D_{\vu}D_{\vv} \Lc \vf\,(\vx)&  \mbox{\textup{ in }} D_1,  \\
      & W (\vx) = 0  &\mbox{ \textup{on} } \partial D_1.
    \end{array}\right.
\end{align*}
%   Note $\operatorname{curl} \vf =\operatorname{curl} \vf^{\,s}$ and $\operatorname{div} \vf^{\,s} =0$, that is,
% \begin{align*}
%   \frac{\partial f^s_2}{\partial x_1}  -  \frac{\partial f^s_1}{\partial x_2} &= \operatorname{curl} \vf  \quad \mbox{ and } \quad      \frac{\partial f^s_1}{\partial x_1} +   \frac{\partial f^s_2}{\partial x_2} = 0.
% \end{align*}
% Using these two relations together with Theorem \ref{th: Inversion of B}, we get
% \begin{align*}
%  \Delta f^s_1 &= -\frac{\partial}{\partial x_2} \operatorname{curl} \vf = -\frac{1}{\det(\vv, \vu)}\frac{\partial}{\partial x_2}D_{\vu}D_{\vv}\, \Lc \vf\\
% \Delta f^s_2 &= \frac{\partial}{\partial x_1} \operatorname{curl} \vf =\frac{1}{\det(\vv, \vu)}\frac{\partial}{\partial x_1}D_{\vu}D_{\vv}\, \Lc \vf.
% \end{align*}
% Therefore, we can recover the Laplacian of each component of $\vf^{\,s}$ from the knowledge of $\Lc \vf$.

%In the next section we prove each of these theorems in two different ways: one of them uses an analytical argument, while the other one presents a geometric explanation.
\section{Proofs of Theorems \ref{th: kernel of B}, \ref{th: kernel of T}, \ref{th: Inversion of B}, \ref{th: Inversion of T}}\label{sec:Proofs}
In this section we prove all four previously stated theorems. We provide two proofs for each one of them: the first proof uses an analytic argument, while the second one presents a geometric explanation.
%\begin{theorem}\label{th: kernel of B}
%The kernel of the V-line Doppler transform $\Bc_{\vu, \vv}$ is set of all  potential vector field $\vf$. In other words,
%\begin{align*}
%\Bc_{\vu,\vv} \vf =0 \quad \mbox{if and only if\ \ } \vf = \nabla V = \left(\frac{\partial V}{\partial x_1}, \frac{\partial V}{\partial x_1}\right) \ \mbox{ for some scalar function } V.
%\end{align*}
%\end{theorem}
\subsection{Proof of Theorem \ref{th: kernel of B}}
%\begin{proof}[Proof of Theorem \ref{th: kernel of B}]
For a given vector field $\vf \in C_c^2(D_1)$, we want to show the existence of a scalar function $V$ satisfying the following:
\begin{align*}
\Lc \vf =0 \quad \mbox{if and only if\ \ } \vf = \nabla V.
\end{align*}
\noindent \textbf{Analytic argument.}\\ Using the definition of $\Lc$ and applying directional derivatives along $\vu$ and $\vv$ we get
\begin{align*}
\Lc \vf= -\mathcal{X}_{\vu}  \left(\vf \cdot \vu\right)  + \mathcal{X}_{\vv}  \left(\vf \cdot \vv\right) &=0\; {\Longrightarrow} \\
   D_{\vu}D_{\vv}\left[-\mathcal{X}_{\vu}  \left(\vf \cdot \vu\right)  + \mathcal{X}_{\vv}  \left(\vf \cdot \vv\right)\right] &=0\;  {\Longrightarrow}\\
  D_{\vv} (\vf \cdot \vu) -D_{\vu} (\vf \cdot \vv) &= 0.
\end{align*}
{The same implications work also in the opposite direction, i.e.
$$  D_{\vv} (\vf \cdot \vu) -D_{\vu} (\vf \cdot \vv) = 0 \quad  \Longrightarrow\quad  \Lc \vf = 0. $$
To see this, consider
\begin{align*}
    D_{\vv} (\vf \cdot \vu) -D_{\vu} (\vf \cdot \vv) &=     D_{\vu}D_{\vv}\left[-\mathcal{X}_{\vu}  \left(\vf \cdot \vu\right)  + \mathcal{X}_{\vv}  \left(\vf \cdot \vv\right)\right]\\
    &=   -D_{\vv}\left\{ D_{\vu}\left[\mathcal{X}_{\vu}  \left(\vf \cdot \vu\right)\right]\right\}  + D_{\vu}\left\{D_{\vv}\left[ \mathcal{X}_{\vv}  \left(\vf \cdot \vv\right)\right]\right\}.
\end{align*}
Since  $\mathcal{X}_{\vu}  \left(\vf \cdot \vu\right)$ and $\mathcal{X}_{\vv}  \left(\vf \cdot \vv\right)$ are constant in the directions $\vu$ and $\vv$ respectively  outside of $D_2$, the functions $D_{\vu}\left[ \mathcal{X}_{\vu}  \left(\vf \cdot \vu\right)\right]$ and $D_{\vv}\left[ \mathcal{X}_{\vv}  \left(\vf \cdot \vv\right)\right]$ are supported in $D_2$. Hence by applying $\Xc_{\vu} \Xc_{\vv}$ to the above equation and using relations \eqref{eq:intertwining relations} , we get
\begin{align*}
    \Xc_{\vu}\Xc_{\vv} \left[   D_{\vv} (\vf \cdot \vu) -D_{\vu} (\vf \cdot \vv) \right] &=     \Xc_{\vu}\Xc_{\vv} \left\{-D_{\vv}\left[ D_{\vu}\left(\mathcal{X}_{\vu}  \left(\vf \cdot \vu\right)\right)\right]  + D_{\vu}\left[D_{\vv}\left( \mathcal{X}_{\vv}  \left(\vf \cdot \vv\right)\right)\right]\right\}  \\
    &= -\mathcal{X}_{\vu}  \left(\vf \cdot \vu\right)  + \mathcal{X}_{\vv}  \left(\vf \cdot \vv\right) \\
    &= \Lc \vf.
\end{align*}
Therefore  $D_{\vv} (\vf \cdot \vu) -D_{\vu} (\vf \cdot \vv) = 0 \quad  \Longrightarrow\quad  \Lc \vf = 0.$ }
% \mj{Note: For a scalar smooth function $D_{\vu}f = 0$ implies $f=0$ if we assume $f$ to be compactly supported. But
% $ -\mathcal{X}_{\vu} \left(\vf \cdot \vu\right)  + \mathcal{X}_{\vv}  \left(\vf \cdot \vv\right) $ is not compactly supported and we need a more detailed argument (the claim still true but I am note sure how short the proof can be.) }

Thus, \begin{align*}
\Lc \vf =0 \quad \mbox{if and only if\ \ } D_{\vv} (\vf \cdot \vu) -D_{\vu} (\vf \cdot \vv) = 0.
\end{align*}
Hence to complete the proof of this theorem it suffices to show that
\begin{align*}
D_{\vv} (\vf \cdot \vu) -D_{\vu} (\vf \cdot \vv) = 0 \quad \mbox{if and only if \ } \vf = \nabla V, \mbox{ for some scalar function } V.
\end{align*}
Consider,
\begin{align}\label{eq:equivalent description for kernel of L}
&D_{\vv} (\vf \cdot \vu) -D_{\vu} (\vf \cdot \vv)\nonumber\\ &= \left(v_1 \frac{\partial}{\partial x_1} +v_2 \frac{\partial}{\partial x_2}\right)(u_1 f_1 +u_2f_2)-\left(u_1 \frac{\partial}{\partial x_1} +u_2 \frac{\partial}{\partial x_2}\right)(v_1 f_1 +v_2f_2)\nonumber \\
&= v_1 u_1\frac{\partial f_1}{\partial x_1} + v_1 u_2 \frac{\partial f_2}{\partial x_1} + v_2 u_1 \frac{\partial f_1}{\partial x_2} +v_2 u_2 \frac{\partial f_2}{\partial x_2}\nonumber
- v_1 u_1\frac{\partial f_1}{\partial x_1} - v_2 u_1 \frac{\partial f_2}{\partial x_1} - v_1 u_2 \frac{\partial f_1}{\partial x_2} - v_2 u_2 \frac{\partial f_2}{\partial x_2}\nonumber \\
&= \det(\vv, \vu) \left(\frac{\partial f_2}{\partial x_1}-\frac{\partial f_1}{\partial x_2}\right)\nonumber\\ &= \det(\vv, \vu)\operatorname{curl}  \vf.
\end{align}
Since $\vu$ and $\vv$ are linearly independent, we conclude
\begin{align*}
D_{\vv} (\vf \cdot \vu) -D_{\vu} (\vf \cdot \vv) = 0 \quad \mbox{if and only if\ \ } \operatorname{curl} \ \vf =0.
\end{align*}
It is known that for {discs}  $\operatorname{curl} \ \vf =0$ if and only if $ \vf = \nabla V$ for some scalar function $V$ {(for instance see \cite[Corollary 16.27] {Lee_Book})}. This completes the proof of Theorem \ref{th: kernel of B}. $\hfill\blacksquare$\\

%%%%%%%%%%%%%%%%%%%%%%%%%%%%%%%%%%%%%%%%%%%%%%

\begin{figure}[ht]
\begin{center}
\begin{subfigure}{.4\textwidth}
\centering
 \includegraphics[height=3.8cm]{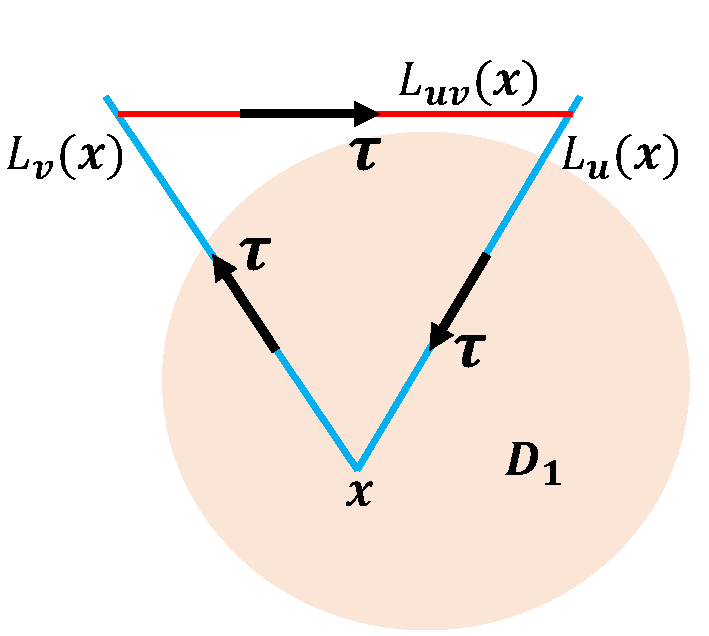}
 \caption{A V-line $L_{\vu}\cup L_{\vv}$ and an additional line segment $L_{\vu\vv}$ outside of $\operatorname{supp} \vf$ with unit tangent vectors $\pmb\tau$.}
  \label{fig2a}
\end{subfigure} \qquad
\begin{subfigure}{.4\textwidth}
\centering
\includegraphics[height=3.8cm]{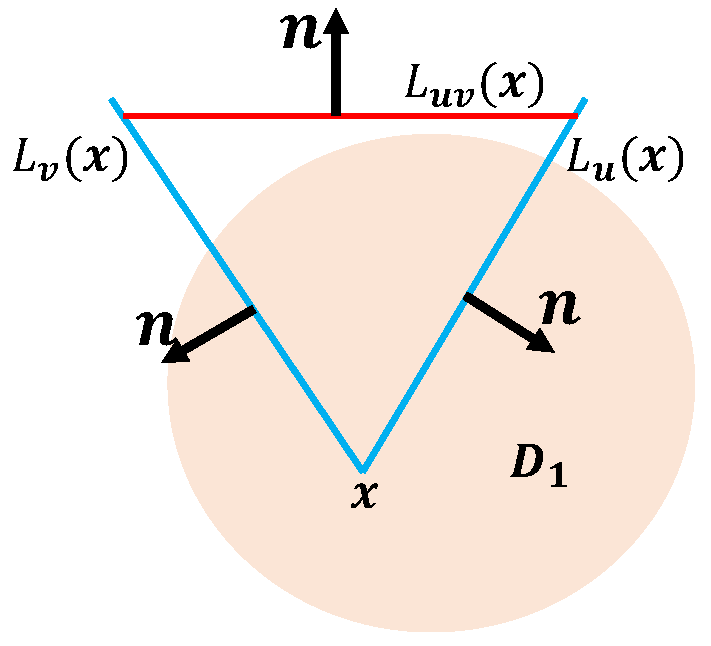}
\caption{A V-line $L_{\vu}\cup L_{\vv}$ and an additional line segment $L_{\vu\vv}$ outside of $\operatorname{supp} \vf$ with unit normal vectors $\vn$.} \label{fig2b}
\end{subfigure}
\end{center}
\vspace{-5mm}
\caption{}
\label{fig2}
\end{figure}

\noindent\textbf{Geometric explanation.}\\
$(\Longleftarrow)$ Assume $\vf = \nabla V$ for some scalar function $V$, thus $\operatorname{curl}\vf=0$.
One can think of transformation $\mathcal{L}$ as the integral of the tangent component of the vector field $\vf$ along branches of the V-lines, i.e.
%------------------------------------------
\begin{align*}
\mathcal{L}\vf = \int_{L_{\vu} \cup L_{\vv}} \vf\cdot \pmb{\tau} \,d t,
\end{align*}
%------------------------------------------
where $\pmb{\tau}$ is the unit tangent vector of the V-line (as shown in Figure \ref{fig2a}).
%----------------------
%%%%%%%%%%%%%%%%%%%%%%%%%%%%%%%%%%%%
%----------------------

Consider a triangular closed  contour defined by some finite intervals of the V-line and an additional ``bridge'' $L_{\vu\vv}$ outside of $D_1\supseteq \operatorname{supp} \vf$ (see Figure \ref{fig2a}). Let $G$ denote the region enclosed by $L_{\vu} \cup L_{\vv} \cup L_{\vu\vv}$.  Using Green's theorem and the fact that $\operatorname{curl}\vf=0$, we get:
\begin{align*}
    \mathcal{L}\vf=\mathcal{L}\vf + \int_{L_{\vu\vv}} \vf \cdot \pmb{\tau} \,d t= \int_{L_{\vu} \cup L_{\vv} \cup L_{\vu\vv}} \vf \cdot \pmb{\tau} \,d t = \int_{G} \operatorname{curl}\vf ds  = 0.
\end{align*}
%sending $L_{\vu\vv}$ segment to infinity and using the decay boundary condition we conclude that $\Lc\vf = 0$.
$(\Longrightarrow)$ The other direction of the statement in Theorem \ref{th: kernel of B} is a consequence of Theorem \ref{th: Inversion of B}. $\hfill\blacksquare$ \\

%\end{proof}
%\begin{theorem}\label{th: kernel of T}
%The kernel of the transverse V-line Doppler transform $\Tc_{\vu, \vv}$ is the set of all vector field $\vf$ that divergence free. In other words,
%\begin{align*}
%\Tc_{\vu,\vv} \vf = 0 \quad \mbox{ if and only if }\ \   \operatorname{div} \vf = \frac{\partial f_1}{\partial x_1}+ \frac{\partial f_2}{\partial x_2}=0.
%\end{align*}
%\end{theorem}
\subsection{Proof of Theorem \ref{th: kernel of T}}
%\begin{proof}[Proof of Theorem \ref{th: kernel of T}]

Recall, we want to prove that a vector field $\vf$ is in the kernel of $\Tc$ if and only if the vector field $\vf$ is divergence-free. \\

\noindent \textbf{Analytic argument.}\\
Due to the following special relation between curl and divergence in $\mathbb{R}^2$:
\begin{equation}\label{eq:curl-div}
\operatorname{curl}\vf^{\perp} = \operatorname{curl}{(-f_2,f_1)} = \frac{\partial f_1}{\partial x_1} - \frac{\partial (-f_2)}{\partial x_2} = \operatorname{div} \vf,
\end{equation}
and the fact that $\mathcal{T}\, \vf = -\mathcal{L}\, \vf^\perp$, Theorem \ref{th: kernel of B} implies
\[
\Lc \vf^{\perp} = 0 \Longleftrightarrow \operatorname{curl}\vf^{\perp} = 0.
\]
Hence
\[
\Tc \vf = 0 \Longleftrightarrow \operatorname{div}\vf = 0,
\]
which ends the proof. $\hfill\blacksquare$\\

%%%%%%%%%%%%%%%%%%%%%%%%%%%%%%%%%%%%%%%%%%%%%%%%

\noindent\textbf{Geometric explanation.}

\noindent$(\Longleftarrow)$ Assume $\operatorname{div} \vf =0$. One can think of $\mathcal{T}\vf$ as the integral of the normal component of the vector field $\vf$ along branches of the V-lines, i.e.
\begin{align*}
\mathcal{T}\vf = \int_{L_{\vu} \cup L_{\vv}} \vf\cdot \vn \,d t,
\end{align*}
%--------------------------------------------
%\begin{figure}[ht]
%\begin{center}
%\begin{subfigure}{.4\textwidth}
%\centering
% \includegraphics[width=3cm]{normal-vline.png}
% \caption{V-lines with unit normal vectors  }
%  \label{fig3a}
%\end{subfigure} \qquad
%\begin{subfigure}{.4\textwidth}
%\centering
%\includegraphics[width=3cm]{normal-vline-closed.png}
%\caption{Attaching a line at infinity to obtain a closed loop} \label{fig3b}
%\end{subfigure}
% \end{center}
%\caption{yyy.}
%\label{fig3}
%\end{figure}
%----------------------------------------------
where $\vn$ is the unit normal vector of the V-line (as shown in Figure \ref{fig2b}).

Consider a triangular closed  contour defined by some finite intervals of the V-line and an additional ``bridge'' $L_{\vu\vv}$ outside of $D_1\supseteq \operatorname{supp} \vf$ (see Figure \ref{fig2b}). Let $G$ denote the region enclosed by $L_{\vu} \cup L_{\vv} \cup L_{\vu\vv}$.  We have
\begin{align*}
    \mathcal{T}\vf=\mathcal{T}\vf + \int_{L_{\vu\vv}} \vf\cdot \vn \,dt = \int_{L_{\vu} \cup L_{\vv} \cup L_{\vu\vv}} \vf \cdot \vn \,dt  = \int_{G} \operatorname{div}\vf ds =0.
\end{align*}
%sending  $L_{\vu\vv}$ segment to infinity and using the decay boundary condition we conclude that $\mathcal{T}\vf = 0$. \\
$(\Longrightarrow)$ The other direction of the statement in Theorem \ref{th: kernel of T} is a consequence of Theorem \ref{th: Inversion of T}. $\hfill\blacksquare$ \\

%%%%%%%%%%%%%%%%%%%%%%%%%%%%%%%%%%%%%%%%%%%%%%%

%\begin{theorem}\label{th: Inversion of B}
%The \operatorname{curl} $\vf$ of a vector field $\vf$ can be recovered uniquely from the knowledge of $\Bc_{\vu,\vv} \vf$.
%\end{theorem}
\subsection{Proof of Theorem \ref{th: Inversion of B}}
%\begin{proof}[Proof of theorem \ref{th: Inversion of B}]
%It is known (e.g. see  \cite[Theorem 3.3.2]{Sharafutdinov_Book}) that every vector field $\vf$ can be decomposed uniquely into a sum of its solenoidal and potential parts as follows:
\noindent \textbf{Analytic argument.}\\
Recall the decomposition of a compactly supported vector field $\vf$ presented in formula (\ref{eq:decomposition}):
\begin{align*}
\vf = \vf^s + \nabla V, \quad \mbox{ with }  \operatorname{div} \vf^s =0\ \mbox{ and } V = 0 \ \mbox{ on }\partial D_1.
\end{align*}
By applying $\Lc$ to this decomposition and using the fact that $\Lc (\nabla V) =0$ (from Theorem \ref{th: kernel of B}), we get
\begin{align*}
%\Lc \vf (\vx) &=  \Lc \vf^s(\vx) + \Lc (\nabla V) (\vx)\\
%\Longrightarrow\qquad \quad
\Lc \vf (\vx) &=  \Lc \vf^s(\vx).
\end{align*}
Taking the directional derivatives $D_{\vu} D_{\vv}$ of the above equation and using formula \eqref{eq:equivalent description for kernel of L} we get
\begin{align}
D_{\vu}D_{\vv} \Lc \vf &= D_{\vu}D_{\vv} \left[-\mathcal{X}_{\vu}  \left(\vf^s \cdot \vu\right)  + \mathcal{X}_{\vv}  \left(\vf^s \cdot \vv\right)\right] \nonumber
=  D_{\vv} (\vf^s \cdot \vu) -D_{\vu} (\vf^s \cdot \vv)\nonumber\\
&= \det(\vv, \vu) \left(\frac{\partial f^s_2}{\partial x_1}- \frac{\partial f^s_1}{\partial x_2}\right). \nonumber
\end{align}
Hence,
\begin{align}
\frac{\partial f^s_2}{\partial x_1}- \frac{\partial f^s_1}{\partial x_2}&= \frac{1}{\det(\vv, \vu)}D_{\vu}D_{\vv} \Lc \vf \label{eq: curl in terms of data}.
\end{align}
Finally, we observe $$\frac{\partial f^s_2}{\partial x_1}- \frac{\partial f^s_1}{\partial x_2} = {\operatorname{curl} \vf^s} = \operatorname{curl} \vf.$$
Combining the last relation with equation \eqref{eq: curl in terms of data} we get the required expression for $\operatorname{curl}\vf$:
\begin{align*}
\operatorname{curl} \vf  &= \frac{1}{\det(\vv, \vu) }D_{\vu}D_{\vv} \Lc \vf.
\end{align*}
This completes the proof of Theorem \ref{th: Inversion of B}. $\hfill\blacksquare$\\
%--------------------------------

\begin{figure}[ht]
\centering
\includegraphics[height=2.6cm]{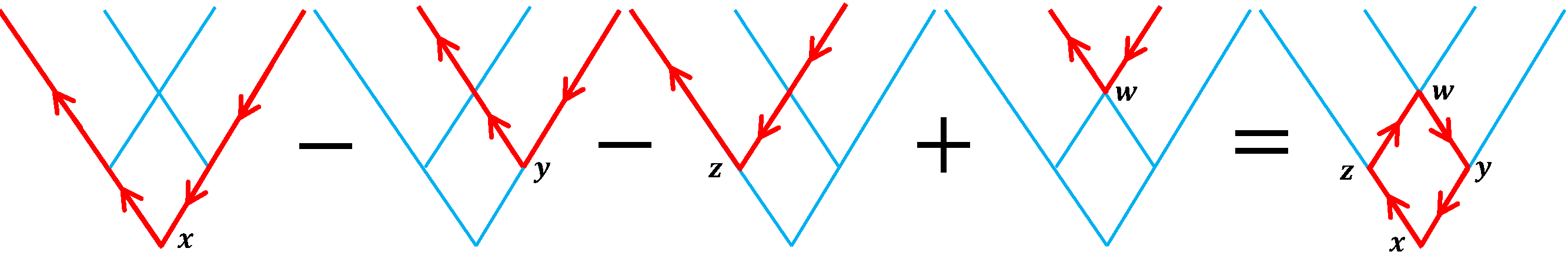}
\caption{A linear combination of $\Lc\vf$ at the vertices of a rhombus resulting in a contour integral of $\vf\cdot\pmb\tau$ along the boundary of the rhombus.}
\label{fig4}
\end{figure}

%%%%%%%%%%%%%%%%%%%%%%%%%%%%%%%%%%%%%%%%%%%%%%%%%%%%%%
\noindent\textbf{Geometric explanation.}\\
Consider the scalar function $h(\vx):=\Lc\vf{(\vx)}$ and the following finite difference of its values at the vertices of a rhombus (refer to Figure \ref{fig4} for visualization):
\begin{equation}
   C\vf{(\vx,\vy,\vz,\vw)}:= \left[h(\vx)-h(\vy)\right]-\left[h(\vz)-h(\vw)\right]=h(\vx)-h(\vy)-h(\vz)+h(\vw).
\end{equation}

\noindent If one sends the side length $\delta>0$ of the rhombus to zero, then
\begin{equation}\label{eq:dir-der}
\lim_{\delta\rightarrow 0}{\left[ \frac{1}{\delta^2}\,C\vf{(\vx,\vy,\vz,\vw)}\right]}= D_{\vu} D_{\vv} h\,(\vx).
\end{equation}

On the other hand, from Figure \ref{fig4} it is easy to see that $C\vf{(\vx,\vy,\vz,\vw)}$ is the clockwise contour integral of $\vf\cdot\pmb\tau$ along the boundary of the rhombus. At the same time, by definition of curl (equivalent to the one in \eqref{eq: definition of div and curl}) for any infinitesimal region $P$ containing $\vx$ we have
\begin{equation}\label{eq:curl-def}
\operatorname{curl} \vf{(\vx)} = \lim_{|P|\rightarrow 0} \frac{1}{|P|} \oint_{\partial P} \vf\cdot \pmb{\tau}\, d t,
\end{equation}
where the integral is taken along the contour traversed counterclockwise.
Since the area of our infinitesimal rhombus is $-\delta^2\det(\vv,\vu)$, formulas \eqref{eq:dir-der} and \eqref{eq:curl-def} imply that
$$
\operatorname{curl} \vf = \frac{1}{\det{(\vv, \vu)}}  D_{\vu} D_{\vv} h,
$$
which is what we wanted to show. $\hfill\blacksquare$\\

%%%%%%%%%%%%%%%%%%%%%%%%%%%%%%%%%%%%%%%%%%%%%%%%%%%%%%%
%%%%%%%%%%%%%%%%%%%%%%%%%%%%%%%%%%%%%%%%
\subsection{Proof of Theorem \ref{th: Inversion of T}}

%\begin{proof}[Proof of theorem \ref{th: Inversion of T}]
\noindent \textbf{Analytic argument.}\\
Using the formula of Theorem \ref{th: Inversion of B} and relation (\ref{eq:curl-div}) between divergence and curl we have
\begin{align}
\operatorname{curl} \vf^\perp &= \frac{1}{\det(\vv, \vu) }D_{\vu}D_{\vv}\, \Lc \vf^\perp, \end{align}
which translates into
\begin{align}
\operatorname{div} \vf &= -\frac{1}{\det(\vv, \vu) }D_{\vu}D_{\vv}\, \Tc \vf,
\end{align}
and concludes the analytic proof of Theorem \ref{th: Inversion of T}. $\hfill\blacksquare$\\
%\end{proof}

%%%%%%%%%%%%%%%%%%%%%%%%%%%%%%%%%%%%%%%%%%%%%%%%%%%%%%%%%%%%%%%%%%%%%%%%%%

\begin{figure}[ht]
\centering
\includegraphics[height=2.6cm]{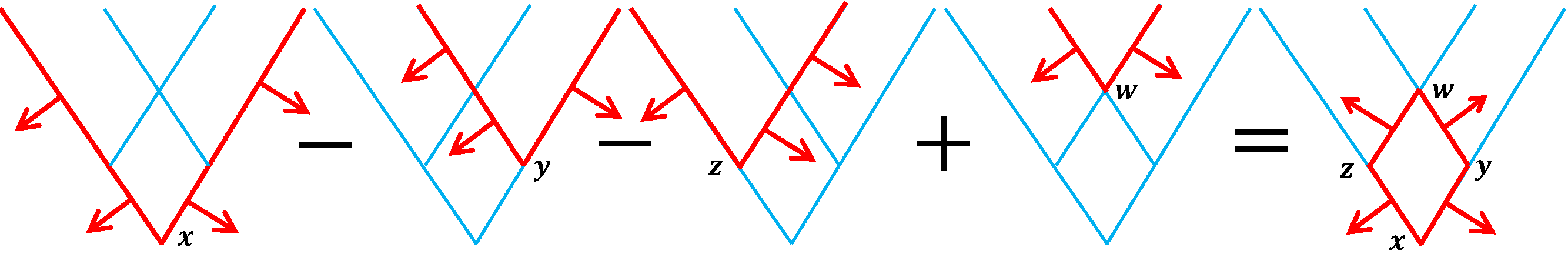}
\caption{A linear combination of $\Tc\vf$ at the vertices of a rhombus resulting in a contour integral of $\vf\cdot \vn$ along the boundary of the rhombus.}
\label{fig5}
\end{figure}

\noindent\textbf{Geometric explanation.}\\
The argument is very similar to that of Theorem \ref{th: Inversion of B}, except $h(\vx):=\Tc\vf{(\vx)}$ here.
Adding and subtracting the values of $h$ as before, we obtain the outward flux of the vector field from the boundary of the infinitesimal rhombus (see Figure \ref{fig5}).

At the same time, by definition of divergence (equivalent to the one in \eqref{eq: definition of div and curl}) for any infinitesimal region $P$ containing $\vx$ we have
$$
\operatorname{div}\vf{(\vx)} = \lim_{|P|\rightarrow 0} \frac{1}{|P|} \oint_{\partial P} \vf\cdot \vn\, d t
$$
Hence,
$$
\operatorname{div} {\vf} = \frac{-1}{\det(\vv, \vu)}  D_{\vu} D_{\vv} h.
$$
Taking $\vf=\nabla V$ completes the proof using $\Delta V = \operatorname{div} {\nabla V}$. $\hfill\blacksquare$

%%%%%%%%%%%%%%%%%%%%%%%%%%%%%%%%%%%%%%%%%%%%%%%%%%%%%%%

\section{Longitudinal and transverse VLT's with their first moments}\label{sec:LVT-moment}
In this section we show that the full vector field $\vf$ can be recovered from the knowledge of its longitudinal V-line transform $\Lc\vf$ and its first moment V-line transform $\Ic\vf$, or alternatively from the knowledge of its transverse V-line transform $\Tc\vf$ and its first moment V-line transform $\Jc\vf$.

The proofs of the theorems presented in this section use the \textbf{signed V-line transform} of a compactly supported scalar function $h$ (see   \cite[Definition 4]{Ambartsoumian_2019}):
\begin{equation}
T_s {h}:=\Xc_{\vu}  h - \Xc_{\vv} h.
\end{equation}
This transform has an explicit inversion formula (e.g. see \cite[Theorem 8]{Ambartsoumian_2019}):
\begin{equation}\label{eq:inv-Ts}
h(\vx)=\frac{1}{||\vv-\vu||}\,D_{\vu}D_{\vv} \int_0^{\infty} (T_sh)(\vx+\vw t)\,dt,
\end{equation}
where
\begin{equation}
   \vw=\frac{\vv-\vu}{||\vv-\vu||}.
\end{equation}

\vspace{3mm}
We can now state and prove the main results of this section.
\begin{thr}\label{th:inversion using moment transform}
Let $\vf$ be a vector field in $\Rb^2$ with components in $C_c^2({D_1})$. Then $\vf$ can be recovered explicitly from  $\Lc \vf$ and $\Ic \vf$.
\end{thr}

\noindent\textbf{Proof.} We know from Theorem \ref{th: Inversion of B} that $\operatorname{curl} \vf$ can be expressed in terms of $\Lc \vf$ as follows:
$$\operatorname{curl} \vf = \frac{1}{\det(\vv, \vu)}\,D_{\vu}D_{\vv}\, \Lc \vf. $$
To prove this theorem we show that the signed V-line transform for each component of $\vf$ can be computed explicitly in terms of $\operatorname{curl} \vf$ and $\Ic \vf$. Indeed,
\begin{align}
\frac{\partial \Ic \vf}{\partial x_1} =  &-\int_0^\infty  t \left( u_1 \frac{\partial f_1}{\partial x_1} +u_2 \frac{\partial f_2}{\partial x_1} \right)(\vx +t \textbf{\textit{u}})\, dt + \int_0^\infty t \left( v_1 \frac{\partial f_1}{\partial x_1} +v_2 \frac{\partial f_2}{\partial x_1} \right)(\vx +t \vv)\, dt\nonumber \\
= &-\int_0^\infty  t \left( u_1 \frac{\partial f_1}{\partial x_1} +u_2 \frac{\partial f_1}{\partial x_2} \right)(\vx +t \textbf{\textit{u}})\, dt -u_2 \int_0^\infty  t \left( \frac{\partial f_2}{\partial x_1} -\frac{\partial f_1}{\partial x_2} \right)(\vx +t \textbf{\textit{u}})\, dt\nonumber\\
&+ \int_0^\infty t \left( v_1 \frac{\partial f_1}{\partial x_1} +v_2 \frac{\partial f_1}{\partial x_2} \right)(\vx +t \vv)\, dt +v_2 \int_0^\infty  t \left( \frac{\partial f_2}{\partial x_1} -\frac{\partial f_1}{\partial x_2} \right)(\vx +t \vv)\, dt\nonumber\\
= &-\int_0^\infty  t\, \frac{d}{dt} f_1(\vx +t \textbf{\textit{u}})\, dt+\int_0^\infty  t\, \frac{d}{dt} f_1(\vx +t \vv)\, dt -u_2\, \Xc^1_{\vu} (\operatorname{curl} \vf) + v_2\, \Xc^1_{\vv} (\operatorname{curl} \vf)\nonumber \\
= & \int_0^\infty  f_1(\vx +t \vu)\, dt - \int_0^\infty f_1(\vx +t \vv)\, dt  -u_2\, \Xc^1_{\vu} {(\operatorname{curl} \vf)} + v_2\, \Xc^1_{\vv} (\operatorname{curl} \vf)\nonumber
\end{align}
In other words,
\begin{equation}\label{eq:relation between I and f1}
\Xc_{\vu} f_1 - \Xc_{\vv} f_1 = \frac{\partial \Ic \vf}{\partial x_1} +u_2\, \Xc^1_{\vu} (\operatorname{curl} \vf) - v_2\, \Xc^1_{\vv} (\operatorname{curl} \vf).
\end{equation}
Differentiating $\Ic \vf$ with respect to $x_2$ and proceeding with a similar calculation, we get
\begin{align}\label{eq:relation between I and f2}
   \Xc_{\vu} f_2 - \Xc_{\vv} f_2 = \frac{\partial \Ic \vf}{\partial x_2} - u_1\, \Xc^1_{\vu} (\operatorname{curl} \vf) + v_1\, \Xc^1_{\vv} (\operatorname{curl} \vf).
\end{align}
Equations \eqref{eq:relation between I and f1} and \eqref{eq:relation between I and f2} express  $T_s{f_1}$ and $T_s{f_2}$ in terms of known $\operatorname{curl} \vf$ and $\Ic \vf$. Therefore, we can recover $f_1$ and $f_2$ explicitly by direct application of formula \eqref{eq:inv-Ts}. $\hfill\blacksquare$
%, which concludes the proof of Theorem \ref{th:inversion using moment transform}.
%$\hfill\blacksquare$

%-------------------------------------------------------

\begin{thr}\label{th:inversion using transverse moment transform}
Let $\vf$ be a vector field in $\Rb^2$ with components in $C_c^2({D_1})$. Then $\vf$ can be recovered explicitly from $\Tc \vf$ and $\Jc \vf$.
\end{thr}
\noindent\textbf{Proof.} From Theorem \ref{th: Inversion of T}, we know that  $\operatorname{div} \vf$ can be expressed in terms of $\Tc \vf$ as follows:
$$\operatorname{div} \vf = -\frac{1}{\det(\vv, \vu)}D_{\vu}D_{\vv} \Tc \vf. $$
In this case, we show that the signed V-line transform of each component of $\vf$ can be computed explicitly in terms of $\operatorname{div} \vf$ and $\Jc \vf$. Indeed, since  $\Jc \vf = - \Ic \vf^{\perp}$ we can use \eqref{eq:relation between I and f1} to get
\begin{align*}
\frac{\partial \Jc \vf}{\partial x_1} = -\frac{\partial \Ic \vf^{\perp}}{\partial x_1}
&=  -\Xc_{\vu} (\vf^\perp)_1 + \Xc_{\vv} (\vf^\perp)_1  +u_2\, \Xc^1_{\vu} (\operatorname{curl} \vf^\perp) - v_2\, \Xc^1_{\vv} (\operatorname{curl} \vf^\perp)\\
&=  \Xc_{\vu} f_2 - \Xc_{\vv} f_2  + u_2\, \Xc^1_{\vu} (\operatorname{div} \vf) - v_2\, \Xc^1_{\vv} (\operatorname{div} \vf),
\end{align*}
where in the last equality we used the relations
\begin{align*}
(\vf^\perp)_1 = -f_2 \mbox{ and } \operatorname{curl} \vf^\perp= \operatorname{div} \vf.
\end{align*}
Therefore, we have
\begin{align*}
\Xc_{\vu} f_2 - \Xc_{\vv} f_2 = \frac{\partial \Jc \vf}{\partial x_1}  - u_2\, \Xc^1_{\vu} (\operatorname{div} \vf) + v_2\, \Xc^1_{\vv} (\operatorname{div} \vf).
\end{align*}
Differentiating $\Jc \vf$ with respect to $x_2$ and proceeding in a similar way, we get
\begin{align*}
\Xc_{\vu} f_1 - \Xc_{\vv} f_1 = -\frac{\partial \Jc \vf}{\partial x_2}  - u_1\, \Xc^1_{\vu} (\operatorname{div} \vf) + v_1\, \Xc^1_{\vv} (\operatorname{div} \vf).
\end{align*}
The last two relations express $T_s{f_1}$ and $T_s{f_2}$ in terms of known $\operatorname{div} \vf$ and $\Jc \vf$. Hence, we can recover $f_1$ and $f_2$ explicitly by direct application of formula \eqref{eq:inv-Ts}. $\hfill\blacksquare$
%%%%%%%%%%%%%%%%%%%%%%%%%%%%%%%%%%%%%%%%%%%%%%%%%%%%%%

\section{Recovery of the full vector field from its star transform}\label{sec:star}
In this section we derive an inversion formula for the star transform of vector-valued functions. Our reconstruction is analogous to the inversion of the star transform of scalar functions introduced in \cite{Amb_Lat_star}.

\begin{defn}\label{def:star_d}
Let $\vgamma_1, \dots, \vgamma_m$ be a distinct  set of unit vectors in $\mathbb{R}^2$. The corresponding \textbf{star transform} $\mathcal{S}\vf$ of a vector %valued function
field
$\vf$ is defined by
\begin{equation}\label{def:star}
\mathcal{S}\vf
%=\sum_{i=1}^m \mathcal{X}_{\vgamma_i}  \begin{bmatrix}
%\langle \vf, \vgamma_i \rangle \\
%\langle \vf, \vgamma_i^{\perp} \rangle
%\end{bmatrix}
=\sum_{i=1}^m c_i \, \mathcal{X}_{\vgamma_i}  \begin{bmatrix}
  \vf \cdot \vgamma_i  \\
  \vf \cdot \vgamma_i^{\perp}
\end{bmatrix}.
\end{equation}
where $c_1, \dots , c_m$ is a set of non-zero weights in $\mathbb{R}$.
\end{defn}

Note that, in contrast with our definition of the V-line transform,  the star transform data contains both the longitudinal and transverse components (this simplifies our discussion). Now, let $\mathcal{R}h(\vpsi,s)$ denote the ordinary Radon transform of a scalar function $h$ in $\mathbb{R}^2$, along the line normal to the unit vector $\vpsi$ and at signed distance $s$ from the origin. Lemmas 1 and 2 in \cite{Amb_Lat_star} provide the following identity:

\begin{equation}\label{lemma12-star}
    \frac{d}{ds}\mathcal{R}(\mathcal{X}_{\vgamma_i}h)(\vpsi,s) = \frac{-1}{ \vpsi\cdot \vgamma_i}\, \mathcal{R}h {(\vpsi, s)}.
\end{equation}

{
\begin{defn} Consider the vector-valued star transform $\mathcal{S}\vf$ with branch directions $\vgamma_1, \dots , \vgamma_m$. We call
\begin{equation}
\mathcal{Z}_1:=\cup_{i=1}^m\left\{\vpsi:\;\vpsi \cdot \vgamma_i= 0\right\}
\end{equation}
the set of singular directions of type 1 for $\mathcal{S}$.

Now let
\begin{equation}
\vgamma(\vpsi) := -\sum_{i=1}^m  \frac{c_i\, \vgamma_i}{ \vpsi \cdot \vgamma_i } \; \in \mathbb{R}^2, \;\;\vpsi\in \mathbb{S}^1\setminus\mathcal{Z}_1.
\end{equation}
We call
\begin{equation}
\mathcal{Z}_2:=\left\{\vpsi:\;\vgamma(\vpsi)= 0\right\}
\end{equation}
the set of singular directions of type 2 for $\mathcal{S}$.
\end{defn}
}

\begin{thr}\label{th:star} {Let $\vf$ be a vector field in $\Rb^2$ with components in $C_c^2(D_1)$}.
Consider {its} vector-valued star transform $\mathcal{S}\vf$ with branch directions $\vgamma_1, \dots , \vgamma_m$ and let

\begin{equation}
Q(\vpsi) := \begin{bmatrix}
 \vgamma(\vpsi) \\
 \vgamma(\vpsi)^{\perp}
\end{bmatrix}^{-1}  \, \in GL(2,\mathbb{R}),\;\;{\vpsi\in \mathbb{S}^1\setminus\left(\mathcal{Z}_1\cup\mathcal{Z}_2\right).}
\end{equation}
%where $\vpsi\in S^1$ is such that $\vgamma(\vpsi)\ne 0$ and $\vpsi \cdot \vgamma_i\ne 0$ for all $i$.
{Then for any $\vpsi\in \mathbb{S}^1\setminus\left(\mathcal{Z}_1\cup\mathcal{Z}_2\right)$ and any $s\in\mathbb{R}$} we have
%\noindent If the unit vector $\vpsi$ is in the domain of $Q(\vpsi)$, then
\begin{equation}\label{star-inv-formula}
Q(\vpsi) \frac{d}{ds}\mathcal{R}(\mathcal{S}\vf) (\vpsi,s) = \mathcal{R}\vf\,{(\vpsi,s)},
\end{equation}
where $\mathcal{R}\vf$ is the component-wise Radon transform of a vector field in $\mathbb{R}^2$.
\end{thr}

\noindent\textbf{Proof}. From (\ref{lemma12-star}) we get

\begin{align}
\frac{d}{ds}\mathcal{R}(\mathcal{S}\vf)  =   \sum_{i=1}^m c_i \frac{d}{ds}\mathcal{R} \mathcal{X}_{\vgamma_i}
\begin{bmatrix}
  \vf \cdot \vgamma_i  \\
  \vf \cdot \vgamma_i^{\perp}
\end{bmatrix} %\nonumber\\
=
-\sum_{i=1}^m  \frac{c_i}{\vpsi \cdot \vgamma_i}\, \mathcal{R}
\begin{bmatrix}
\vf \cdot \vgamma_i \\
\vf \cdot \vgamma_i^{\perp}
\end{bmatrix}
\end{align}
Using the linearity of $\mathcal{R}$ and inner product, we simplify the last expression further to obtain
\begin{equation}\label{finally}
\frac{d}{ds}\mathcal{R}(\mathcal{S}\vf) =  \begin{bmatrix}
 \mathcal{R}\vf \cdot \vgamma(\vpsi) \\
 \mathcal{R}\vf \cdot \vgamma(\vpsi)^{\perp}
\end{bmatrix} = Q(\vpsi)^{-1} \mathcal{R}\vf.
\end{equation}
Finally,
\begin{equation*}
Q(\vpsi) \frac{d}{ds}\mathcal{R}(\mathcal{S}\vf) (\vpsi,s) = \mathcal{R}\vf(\vpsi,s),
\end{equation*}
which ends the proof.$\hfill\blacksquare$\\

{It is easy to see that if $\mathcal{Z}_1\cup\mathcal{Z}_2$ is finite, then all singularities appearing in the left hand side of formula \eqref{star-inv-formula} are removable. In other words,
\begin{rem}
 If $\left\{\vzeta^{(i)}\right\}_{i=1}^M=\mathcal{Z}_1\cup\mathcal{Z}_2$, $M\ge m$,
then by formula \eqref{finally} and continuity of $\mathcal{R}\vf$ we have

\begin{equation}
 \lim_{\vpsi\to\vzeta^{(i)}}\left[Q(\vpsi) \frac{d}{ds}\mathcal{R}(\mathcal{S}\vf)(\vpsi,s)\right] =\lim_{\vpsi\to\vzeta^{(i)}}\left[Q(\vpsi) Q(\vpsi)^{-1} \mathcal{R}\vf\, (\vpsi,s) \right]  = \mathcal{R}\vf\,\left(\vzeta^{(i)},s\right)<\infty.
\end{equation}
%Hence, one can apply $\mathcal{R}^{-1}$ to both sides of equation \eqref{star-inv-formula} to recover $\vf\,{(\vx)}$ at every $\vx\in D_1$.
\end{rem}

\begin{cor}
If $\mathcal{Z}_1\cup\mathcal{Z}_2$ is finite, then one can apply $\mathcal{R}^{-1}$ to both sides of equation \eqref{star-inv-formula} to recover $\vf\,{(\vx)}$ at every $\vx\in D_1$.
\end{cor}
}

{It is clear that the set $\mathcal{Z}_1$ is finite.}
\begin{comment}
Now, if we assume that $ \vpsi \cdot \vgamma_i \neq 0$,  it is easy to notice that matrix function $Q(\vpsi)$ is undefined if and only if
\begin{equation}\label{eq:pol-zero}
\vgamma(\vpsi) = -\sum_{i=1}^m  \frac{c_i\, \vgamma_i}{ \vpsi \cdot \vgamma_i}=0.
\end{equation}
\end{comment}
{
Let us study in more detail the set $\mathcal{Z}_2$, i.e. the solutions of the equation
\begin{equation}\label{eq:pol-zero}
\vgamma(\vpsi) = -\sum_{i=1}^m  \frac{c_i\, \vgamma_i}{ \vpsi \cdot \vgamma_i}=0.
\end{equation}}
{Bringing the fractions in the above sum to a common denominator, we have

\begin{equation}\label{psiDenum}
      \vgamma(\vpsi) =  \frac{ -\sum_{i=1}^m \left(  \prod_{j \neq i} \vpsi \cdot \vgamma_j \right) c_i \vgamma_i }{\prod_{j = 1}^m \vpsi \cdot \vgamma_j}.
\end{equation}
The numerator is a vector function, which we denote by $-\vP(\vpsi)$.

Using the notation $\vgamma_j=(a_j,b_j)$, $j=1,\ldots,m$, we can re-write $\vP(\vpsi)$ in terms of components of $\vpsi$ as follows
\begin{equation}\label{P2rs1}
      \vP(\psi_1,\psi_2)=\sum_{i=1}^m \left[c_i \prod_{j \neq i} (\psi_1 a_j + \psi_2 b_j)\right] (a_i,b_i).
\end{equation}
Notice, that each component {$P_i$ of  $\vP=(P_1,P_2)$} is a homogeneous polynomial in terms of $\vpsi=(\psi_1,\psi_2)$. It is not hard to see {that on $\mathbb{S}^1\setminus \mathcal{Z}_1$} the condition $\vgamma(\vpsi)=0$ reduces to $\vP(\vpsi)=0$, {i.e. $P_1(\vpsi)=0$ and $P_2(\vpsi)=0$}.
}

{The classical B\'ezout's Theorem from algebraic geometry implies that: two affine algebraic plane curves $C_1$ and $C_2$ of degree $d_1$ and $d_2$ respectively either have a common component (i.e., the polynomials defining them have a common factor) or intersect in at most $d_1d_2$ points (e.g. see \cite{Fischer}).
Since $P_c(\psi_1,\psi_2):=\psi_1^2+\psi_2^2-1$ is irreducible, each homogeneous polynomial $P_i(\psi_1,\psi_2)$ either has finitely many zeros on $\mathbb{S}^1$, or is identically zero. In other words, the set $\mathcal{Z}_2$ is either finite or contains $\mathbb{S}^1\setminus\mathcal{Z}_1$. Thus,
}

\begin{cor}\label{cor1}
The star transform is invertible if  $\vgamma(\vpsi)\ne 0$ for at least one {$\vpsi\in \mathbb{S}^1\setminus\mathcal{Z}_1$. In fact, the proof of Theorem \ref{th:injectivity} will show, that the converse implication is also true, i.e. it is an ``if and only if'' statement.}
\end{cor}
This corollary help us to give a complete description of invertible star configurations.
\begin{defn}
 We call a star transform { $\mathcal{S}$} \textbf{symmetric}, if  $m=2k$ for some $k\in\mathbb{N}$ and (after possible re-indexing) $\vgamma_i=-\vgamma_{k+i}$  with $c_i=-c_{k+i}$ for all $i=1,\ldots, k$.
\end{defn}

As a side note, the sign convention in $c_i=-c_{k+i}$ is different from that of \cite{Amb_Lat_star}, which is due to the orientation that we are using in the definition of the star transform for vector fields.

\begin{thr}\label{th:injectivity}
The star transform { $\mathcal{S}$} is  invertible if and only if it is not symmetric.
\end{thr}
\noindent \textbf{Proof.} The argument follows closely the steps of the proof of Theorem 2 in \cite{Amb_Lat_star}, which we present here for completeness.
{If $\mathcal{S}$ is symmetric with $2k$ ray directions, then the star transform data {contains less information than the data of standard Radon transform in $k$ fixed directions (one can obtain the star transform from a sum of corresponding pairs of the Radon transform)}. It is well known that we cannot recover a function from Radon data of finitely many angles, therefore, $\mathcal{S}$ cannot be invertible.}

Now, assume that for a fixed choice of  $\vgamma_1, \dots, \vgamma_m$ there is no inversion for the corresponding star transform. By Corollary \ref{cor1} we have $\vgamma(\vpsi)\equiv 0$. Hence, the $\vP(\psi_1,\psi_2)\equiv0$.

   { Without loss of generality assume that $a_1 \neq 0$ (otherwise $b_1\neq 0$).} If we take $\psi_1=b_1$ and $\psi_2=-a_1$, then formula \eqref{P2rs1} implies the following
    \[
    %0 = e_{m-1}(ra+sb)
    0= {a_1} \prod_{j \neq 1} (\psi_1 a_j + \psi_2 b_j) = {a_1}(a_2 b_1 - a_1 b_2) \dots (a_m b_1 - a_1 b_m).
    \]
       Hence, for some index {$\ell$} we are required to have {$a_{\ell} b_1 = a_1 b_{\ell}$} or equivalently { $\frac{a_1}{b_1} = \frac{a_{\ell}}{b_{\ell}}$ }. Given the assumption that $\vgamma_i$'s are distinct unit vectors, we conclude that {$\vgamma_1 = - \vgamma_{\ell}$}.  Applying this procedure with $\psi_1=a_j$ and $\psi_2=-a_j$ repeatedly for all $j$, we conclude that in order for the star transform to be non-invertible,  its ray directions have to come in opposite pairs. \\

Now, we prove the relation $c_{i}=-c_{k+i}$ between the corresponding weights of each pair.  Without loss of generality, let us assume that $m=2k$, $\vgamma_i=-\vgamma_{k+i}$ for $i=1,\ldots,k$ and no pair of vectors $\vgamma_1,\ldots,\vgamma_k$ are collinear. Then we can re-write formula (\ref{P2rs1}) as
\begin{equation*}%\label{P2rs2}
      \vP(\psi_1,\psi_2)=(-1)^{k}\prod_{j=1}^k (\psi_1 a_j + \psi_2 b_j)\sum_{i=1}^k \left[(c_i+c_{k+i})\,\vgamma_i \prod_{\substack{j=1 \\ j\neq i}}^k (\psi_1 a_j + \psi_2 b_j)\right] = 0, \; \forall \,\; \psi_1,\psi_2\in \mathbb{R}.
\end{equation*}
Since $\prod_{j=1}^k (\psi_1 a_j + \psi_2 b_j)= (\vpsi\cdot \vgamma_1)  \dots  (\vpsi \cdot \vgamma_k)$, we have
\begin{equation}\label{q2-noninj}
 \vgamma(\vpsi)= {\displaystyle\frac{-\vP(\psi_1,\psi_2)}{(\vpsi\cdot \vgamma_1) \dots (\vpsi \cdot \vgamma_m)}}=\frac{-1}{(\vpsi\cdot \vgamma_1) \dots (\vpsi \cdot \vgamma_k)}\;{\displaystyle \sum_{i=1}^k \left[(c_i+c_{k+i})\,\vgamma_i \prod_{\substack{j=1 \\ j\neq i}}^k (\psi_1 a_j + \psi_2 b_j)\right]}
\end{equation}
for all $\vpsi\in S^1$ outside of the finite set $\left\{\vpsi:\vpsi\cdot \vgamma_i = 0,\;i=1,\ldots,m\right\}$.

Hence, in order for $\mathcal{S}$ to be non-invertible, we must have
\begin{equation*}
   \sum_{i=1}^k \left[(c_i+c_{k+i})\,\vgamma_i \prod_{\substack{j=1 \\ j\neq i}}^k (\psi_1 a_j + \psi_2 b_j)\right] \equiv0.
\end{equation*}
Following the argument from the first part of the proof, if we take $\psi_1=b_1$ and $\psi_2=-a_1$, then
    \[
    %0 = e_{m-1}(ra+sb)
    0= \left(c_1+c_{k+1}\right)\prod_{j=2}^k(\psi_2 a_j + \psi_2 b_j)  = \left(c_1+c_{k+1}\right)(a_2 b_1 - a_1 b_2) \dots (a_k b_1 - a_1 b_k).
    \]
    Since all vectors $\vgamma_1,\ldots, \vgamma_k$ are pairwise linearly independent, $a_jb_1-a_1b_j\ne0$ for $j=2,\ldots,k$. Hence, the last equation implies $c_{k+1}=-c_1$.

    Applying this procedure with $\psi_1=a_j$ and $\psi_2=-a_j$ repeatedly for all $j$, we conclude that in order for the star transform to be non-invertible, we must have $c_{k+i}=-c_i$ for all $i=1,\ldots, k$.
 $\hfill\blacksquare$

 Theorem \ref{th:injectivity} immediately implies the following.
%---------------------------------------------
\begin{cor}
 Any vector-valued star transform with an odd number of rays is invertible.
\end{cor}

{
\begin{rem}
It is easy to notice, that similar to the case of V-line transforms, there exists a disc $D_2$ of finite radius $r_2>0$ with the property that the restriction of $\Sc\vf$ to $\overline{D_2}$ completely defines it in $\mathbb{R}^2$. In other words, to recover a vector field $\vf$ in $\Rb^2$ with components in $C_c^2(D_1)$ one only needs to know its vector-valued star transform $\Sc\vf$ in $D_2$.
\end{rem}
}

%\vspace{1cm}
%Let us rewrite equation \eqref{eq:pol-zero} using complex exponentials and interpreting $\mathbb{R}^2$ as $\mathbb{C}$. In other words, we rewrite $\vpsi=e^{i\beta}$ and $\vgamma_k=e^{i\alpha_k}$ to get:
%\begin{equation}\label{eq:complex-plane}
%\vgamma(\vpsi) = -\sum_{k=1}^m  \frac{c_k\, %e^{i\alpha_k}}{ \cos{(\alpha_k-\beta)}}=-e^{i\beta} \sum_{k=1}^m  \frac{c_k\, e^{i(\alpha_k-\beta)}}{ \cos{(\alpha_k-\beta)}}=0.
%\end{equation}
%Let $\vpsi_0$ and the corresponding polar angle $\beta_0$ be such that $\cos{(\alpha_k-\beta_0)}\ne0$ for all $k=1,\ldots,m$. Then, from equation \eqref{eq:complex-plane} and Corollary \ref{cor1} we get
%\begin{cor}\label{cor2}
%The star transform $\mathcal{S}\vf$ is invertible whenever $c_1+\ldots+c_m \ne0$.
%\end{cor}

\begin{rem}
Similar to the case of the star transform of scalar functions studied in \cite{Amb_Lat_star}, in invertible configurations the function  $Q(\vpsi)$ (or equivalently the function $\vgamma(\vpsi)$) contains information about stability of inversion of the  star transform on vector fields. A comprehensive analysis of that function is a non-trivial task (e.g. see \cite{Amb_Lat_star} for a similar problem in the scalar case) and the authors plan to address it in another publication.
\end{rem}

\begin{rem}
When $m=2$ and $c_1=-c_2=1$, the star transform of vector field $\vf$ corresponds to the vector function $(\Lc\vf,\Tc\vf)$. Hence, Theorem \ref{th:star} provides another approach to recovering the full vector field $\vf$ from its longitudinal and transverse V-line transforms. In the special case when $\vgamma_1 = -\vgamma_2$ (and only in that case), the matrix $Q(\vpsi)$ is undefined for any $\vpsi$ and the corresponding transform is not invertible.
\end{rem}

%%%%%%%%%%%%%%%%%%%%%%%%%%%%%%%%%%%%%%%%%%%%%%%%%%%

\section{Additional remarks}\label{sec:remarks}
\begin{enumerate}
\item There are some interesting similarities between the V-line vector tomography and classical (straight line) vector tomography, despite the differences in the concepts and techniques of deriving the results.
%\textbf{Similarities:}
\begin{itemize}
\item The kernel descriptions for the longitudinal and transverse transforms are identical in the V-line case (obtained in this paper) and the straight line case (see \cite{Derevtsov2}).
\item In article \cite{Derevtsov2}, the authors showed that the combination of LRT and TRT provides a unique reconstruction of a vector field in $\Rb^2$. We achieve the same result combining the V-line versions of those transforms.
\item The authors in \cite{Venky_Suman_Manna, R_K_Mishra} used the  combination of LRT and the first integral moment transform data to get the full vector field in $\Rb^2$. Our Theorem \ref{th:inversion using moment transform} achieves the same result in the case of V-line transforms.
\end{itemize}

\item In dimensions $n \geq 3$, the longitudinal $V$-line transform can be defined in the same fashion as for $n=2$, but the transverse $V$-line transform will require more details, since there is an $(n-1)$-dimensional space of transverse directions. Once a proper choice for the transverse direction is made, techniques similar to the ones introduced in this paper can be used to study injectivity and invertibility for both transforms  in higher dimensions. The authors plan to address these questions in a future publication.

{
\item This article contains multiple inversion formulas for various integral transforms of a compactly supported vector field in $\mathbb{R}^2$. Numerical implementation of these formulas and corresponding algorithms, as well as the study of their stability and artifacts are challenging tasks (e.g. see \cite{Amb_Lat_star} for discussion of stability issues of inverting the (scalar) star transform). The authors plan to address these questions in a future publication.

\item A typical application of vector field tomography is the recovery of fluid flow from reciprocal measurements of certain signals through the domain of the flow. Some examples of such measurements include changes in travel time of acoustic signals sent in opposite directions or changes in the optical path length (integral of the refraction index along the linear path) of a laser beam shined through the fluid \cite{Norton}. The latter model assumes that the photon beams are collimated and the changes of the (linear) optical path are expressed in the optical phase shift measured interferometrically. In this setup the information carried by reflected or scattered photons is lost. At the same time, if the fluid has small suspended particles with a different refractive index than the medium, they will cause some photons change their direction due to reflection and scattering. It is conceivable that one can use a sensor array on the opposite (from the laser source) side of the flow region to measure the data carried by such reflected particles. In the case of single reflections, a potential model for such data can be based on V-line transforms of the flow field.
}

%\item Many of the results and techniques of this paper can be generalized naturally to a large class of Riemannian surfaces with well defined V-line transforms. In this  general setting, we fix two branch directions $\vu$ and $\vv$ on the surface using a connection. The transformations are defined by integrating over the geodesics in these directions starting from any given point. The parallel vector fields defined by $\vu$ and $\vv$ will then play the rule of the directional derivatives $D_{\vu}$ and $D_{\vv}$ appearing in the inversion formulas.
\end{enumerate}

%------------------------------------------

\section{Acknowledgements}\label{sec:acknowledge}
The work of G. Ambartsoumian was partially funded by NSF grant DMS 1616564. {He would also like to thank Dr. Matthew Lewis for useful discussions about potential applications of this work.} R. K. Mishra would like to thank Dr. Souvik Roy for providing postdoctoral funding at the University of Texas at Arlington.

%------------------------------------------
%\vspace{5mm}

%\bibliographystyle{plain}
%\bibliography{reference}

\end{document}

%% file: preamble.tex
\newcommand{\cout}[1]{}

%% file: BDT.bbl
\begin{thebibliography}{100}

\bibitem{Anuj_TRT}
Anuj Abhishek.
\newblock Support theorems for the transverse ray transform of tensor fields of
  rank $m$.
\newblock {\em Journal of Mathematical Analysis and Applications}, 485(2),
  2020.

\bibitem{Anuj_Rohit}
Anuj Abhishek and Rohit~Kumar Mishra.
\newblock Support theorems and an injectivity result for integral moments of a
  symmetric $m$-tensor field.
\newblock {\em Journal of Fourier Analysis and Applications}, 25(4):1487--1512,
  August 2019.

\bibitem{Gaik_2012}
Gaik Ambartsoumian.
\newblock Inversion of the V-line Radon transform in a disc and its  applications in imaging.
\newblock {\em Computers $\&$ Mathematics with Applications}, 64(3):260 -- 265,  2012.
\newblock Mathematical Methods and Models in Biosciences.

\bibitem{Amb2} Gaik Ambartsoumian.
\newblock V-line and conical Radon transforms with applications in imaging,
\newblock Chapter 7 in \textit{``The Radon Transform: The First 100 Years and Beyond''}, edited by R. Ramlau and O. Scherzer, Radon Series on Computational and Applied Mathematics, De Gruyter, 2019.

\bibitem{Ambartsoumian_2019}
Gaik Ambartsoumian and Mohammad J. Latifi Jebelli.
\newblock The V-line transform with some generalizations and cone  differentiation.
\newblock {\em Inverse Problems}, 35(3):034003, Feb 2019.

\bibitem{Amb_Lat_star}
Gaik Ambartsoumian and Mohammad J. Latifi Jebelli.
\newblock Inversion and symmetries of the star transform,
\newblock {\em arXiv} preprint arXiv:2005.01918

\bibitem{Gaik_Moon}
Gaik Ambartsoumian and Sunghwan Moon.
\newblock A series formula for inversion of the V-line {R}adon transform in a
  disc.
\newblock {\em Computers $\&$ Mathematics with Applications}, 66(9):1567 --
  1572, 2013.
\newblock BioMath 2012.

\bibitem{Amb-Roy}
Gaik Ambartsoumian and Souvik Roy.
\newblock Numerical inversion of a broken ray transform arising in single scattering optical tomography,
\newblock \textit{IEEE Transactions on Computational Imaging},  2(2)  (2016), pp 166--173.

\bibitem{Denisjuk1994}
Alexander Denisjuk.
\newblock Inversion of generalized Radon transform.
\newblock {\em Am. Math. Soc. Trans}, 1994.

\bibitem{Denisjuk_Paper}
Alexander Denisjuk.
\newblock Inversion of the X-ray transform for 3{D} symmetric tensor fields  with sources on a curve.
\newblock {\em Inverse Problems}, 22(2):399--411, 2006.

\bibitem{Derevtsov2}
Evgeny~Yu.~Derevtsov and Valery~V.~Pickalov.
 Reconstruction of vector fields and their singularities from ray transform. {\em Numerical Analysis and Applications}, 4(1): 21--35, 2011.

\bibitem{Derevtsov3}
Evgeny~Yu.~Derevtsov and Ivan~E.~Svetov. Tomography of tensor fields in the plane.
{\em Eurasian J. Math. Comp. Applications}, 3(2):24--68, 2015.

\bibitem{Fischer}
Gerd Fischer.
\newblock {\em Plane Algebraic Curves}.
\newblock American Mathematical Society, 2001.


\bibitem{Florescu_Markel_Schotland_2011}
Lucia Florescu, Vadim~A. Markel and John~C. Schotland.
\newblock Inversion formulas for the broken-ray Radon transform.
\newblock {\em Inverse Problems}, 27(02):025002, 2011.

\bibitem{Florescu_Markel_Schotland_2009}
Lucia Florescu, John~C. Schotland and Vadim~A. Markel.
\newblock Single-scattering optical tomography.
\newblock {\em Phys. Rev. E}, 79(3):036607, 2009.

\bibitem{Florescu3}
Lucia Florescu, John~C. Schotland and Vadim~A. Markel. \newblock Single scattering optical tomography: simultaneous reconstruction of scattering and absorption. \newblock {\em Phys. Rev. E},  81(1):016602, 2010.

\bibitem{Florescu4} Lucia Florescu, John~C. Schotland and Vadim~A. Markel.
\newblock Nonreciprocal broken ray transforms with applications to fluorescence imaging.
\newblock {\em Inverse Problems}  34(9):094002, 2018.

\bibitem{Rim_Gaik_2014}
Rim Gouia-Zarrad and Gaik Ambartsoumian.
\newblock Exact inversion of the conical Radon transform with a fixed opening
  angle.
\newblock {\em Inverse Problems}, 30(4):045007, March 2014.

\bibitem{Rohit_Griemaier}
Roland Griesmaier, Rohit~Kumar Mishra and Christian Schmiedecke.
\newblock Inverse source problems for {M}axwell's equations and the windowed
  {F}ourier transform.
\newblock {\em SIAM J. Sci. Comput.}, 40(2):A1204--A1223, 2018.

%\bibitem{Haltmeier_Moon_Schiefeneder}
%Markus Haltmeier, Sunghwan Moon and Daniela Schiefeneder.
%\newblock Inversion of the attenuated V-line transform with vertices on the circle.
%\newblock {\em IEEE Transactions on Computational %Imaging}, 3(4):853--863, December 2017.

%\bibitem{Haltmeier_2014}
%Markus Haltmeier.
%\newblock Exact reconstruction formulas for a Radon transform over cones.
%\newblock {\em Inverse Problems}, 30(3):035001, Feb 2014.


\bibitem{Holman2013}
Sean Holman.
\newblock Generic local uniqueness and stability in polarization tomography.
\newblock {\em Journal of Geometric Analysis}, 23:229--269, 2013.



\bibitem{Hoop2019}
Maarten V. de Hoop, Teemu Saksala and Jian Zhai.
\newblock Mixed ray transform on simple 2-dimensional Riemannian manifolds.
\newblock {\em Proceedings of the American Mathematical Society.} 147, 4901-4913, 2019.
%\bibitem{Hristova_2015}
%Yulia Hristova.
%\newblock Inversion of a V-line transform arising in emission tomography.
%\newblock {\em Journal of Coupled Systems and Multiscale Dynamics,},  3(3):272--277, September 2015.

%\begin{comment}%%%%%%%%%%%%%%%%%%%%%%%%%%%%%%%%%
%\bibitem{Hubenthal}
%Mark Hubenthal.
%\newblock The broken ray transform on the square.
%\newblock {\em Journal of Fourier Analysis and Applications}, 20:1050--1082,2014.

%\bibitem{Hubenthal_2015}
%Mark Hubenthal.
%\newblock The broken ray transform in $n$-dimensions with flat reflecting
%  boundary.
%\newblock {\em Inverse Problems $\&$ Imaging}, 9(1):143--161, 2015.

%\bibitem{Ilmavirta_2013}
%Joonas Ilmavirta.
%\newblock Broken ray tomography in the disc.
%\newblock {\em Inverse Problems}, 29(3):035008, Mar 2013.

%\bibitem{Ilmavirta_Salo_2016}
%Joonas Ilmavirta and Mikko Salo.
%\newblock Broken ray transform on a Riemann surface with a convex obstacle.
%\newblock {\em Communications in Analysis and Geometry}, 24(02):379 – 408,
%  June 2016.
%\end{comment}%%%%%%%%%%%%%%%%%%%%%%%%%%%%%%%%

%\bibitem{Jung_2015}
%Chang-Yeol Jung and Sunghwan Moon.
%\newblock Inversion formulas for cone transforms arising in application of  Compton cameras.
%\newblock {\em Inverse Problems}, 31(1):015006, Jan 2015.

%\bibitem{Jung_Moon_yeol}
%Chang-Yeol Jung and Sunghwan Moon.
%\newblock Exact inversion of the cone transform arising in an application of a Compton camera consisting of line detectors.
%\newblock {\em SIAM Journal on Imaging Sciences}, 9(2):520--536, 2016.

\bibitem{Katsevich_2013}
Alexander Katsevich and Roman Krylov.
\newblock Broken ray transform: inversion and a range condition.
\newblock {\em Inverse Problems}, 29(7):075008, June 2013.

\bibitem{Katsevich2006}
Alexander Katsevich.
\newblock Improved cone beam local tomography.
\newblock {\em Inverse Problems}, 22(2):627, 2006.

\bibitem{Katsevich2013}
Alexander Katsevich and Thomas Schuster.
\newblock An exact inversion formula for cone beam vector tomography.
\newblock {\em Inverse Problems}, 29(6):065013, 2013.

\bibitem{Norton} Stephen J. Norton.
\newblock Unique tomographic reconstruction of vector fields using boundary data.
\newblock {\em IEEE Transactions on Image Processing}, 1(3):406--412, July 1992.


\bibitem{Wongsason_2020}
Dojin Kim  and Patcharee Wongsason.
\newblock Three-dimensional vector field inversion formula using first moment transverse transform in quaternionic approaches.
\newblock \textit{Mathematical Methods in the Applied Sciences}, 2020.
\newblock https://doi.org/10.1002/mma.6427

\bibitem{Venky_Suman_Manna}
Venkateswaran~P. Krishnan, Ramesh Manna, Suman~Kumar Sahoo and Vladimir~A.  Sharafutdinov.
\newblock Momentum ray transforms.
\newblock {\em Inverse Problems $\&$ Imaging}, 13(3):679--701, 2019.

\bibitem{Venky_Suman_Manna2}
Venkateswaran~P. Krishnan, Ramesh Manna, Suman~Kumar Sahoo and Vladimir~A.  Sharafutdinov.
\newblock Momentum ray transforms, II: range characterization in the Schwartz
  space.
\newblock {\em Inverse Problems}, 36(4), 2020.

\bibitem{Venky_and_Rohit}
Venkateswaran~P. Krishnan and Rohit~K. Mishra.
\newblock Microlocal analysis of a restricted ray transform on symmetric
  {$m$}-tensor fields in {$\mathbb{R}^n$}.
\newblock {\em SIAM J. Math. Anal.}, 50(6):6230--6254, 2018.

\bibitem{Francois_Rohit_Venky}
Venkateswaran~P. Krishnan, Rohit~K. Mishra, and Fran\c{c}ois Monard.
\newblock On solenoidal-injective and injective ray transforms of tensor fields
  on surfaces.
\newblock {\em Journal of Inverse and Ill-posed Problems}, 27(4):527--538, September 2019.

\bibitem{VRS}
Venkateswaran~P. Krishnan, Rohit~K. Mishra and Suman~K. Sahoo.
\newblock Microlocal inversion of a 3-dimensional restricted transverse ray  transform of symmetric $m$-tensor fields.
\newblock {\em arXiv} preprint arXiv:1904.02812.

\bibitem{Kats-Kryl-15}
Roman Krylov and Alexander Katsevich,
\newblock Inversion of the broken ray transform in the case of energy dependent attenuation,
\newblock\textit{Physics in Medicine \& Biology}, 60(11):4313--4334, 2015.

\bibitem{Lee_Book}
John M. Lee.
\newblock {\em Introduction to
smooth manifolds}, second edition.
\newblock {Graduate Texts in Mathematics, Springer, New York.}



% \bibitem{Schuster_Derevtsov}
% Alfred K.~Louis, Evgeny Y.~Derevtsov, Anton V.~Efimov and Thomas Schuster.
% \newblock Singular value decomposition and its application to numerical
%   inversion for ray transforms in 2D vector tomography.
% \newblock {\em Journal of Inverse and Ill-posed Problems}, 19:689 – 715, 2011.

\bibitem{R_K_Mishra}
Rohit K. Mishra.
\newblock Full reconstruction of a vector field from restricted Doppler and  first integral moment transforms in $\mathbb{R}^n$.
\newblock {\em Journal of Inverse and Ill-posed Problems}, 2019.

\bibitem{Rohit_Suman_2020}
Rohit K. Mishra and Suman K. Sahoo.
\newblock Injectivity and range description of first $(k+1)$ integral moment transforms over $m$-tensor fields in $\mathbb{R}^n$.
\newblock {\em Preprint,} Arxiv:2006.13102, 2020.

\bibitem{Monard2016a}
François Monard.
\newblock Efficient tensor tomography in fan-beam coordinates.
\newblock {\em Inverse Problems $\&$ Imaging}, 10:433--459, 2016.

\bibitem{Nattrer_Wuebbeling}
Frank Natterer and Frank Wübbeling.
\newblock {\em Mathematical Methods in Image Reconstruction}.
\newblock Society for Industrial and Applied Mathematics, 2001.

\bibitem{Novikov_Sharafutdinov}
Roman Novikov and Vladimir Sharafutdinov.
\newblock On the problem of polarization tomography: I.
\newblock {\em Inverse Problems}, 23(3), 2007.

\bibitem{Palamodov2009}
Victor Palamodov.
\newblock Reconstruction of a differential form from Doppler transform.
\newblock {\em SIAM Journal on Mathematical Analysis}, 41(4):1713--1720, 2009.

\bibitem{PSU2012}
Gabriel P. Paternain, Mikko Salo and Gunther Uhlmann.
\newblock Tensor tomography on surfaces.
\newblock {\em Inventiones Mathematicae} 193, 229-247, 2013.


\begin{comment}
\bibitem{Bezout}
Joachim Schmid.
\newblock On the affine Bezout inequality.
\newblock {\em Manuscripta Math} 88, pages 225-232, 1995.
\end{comment}

\bibitem{Schuster2000}
Thomas Schuster.
\newblock The 3D Doppler transform: elementary properties and computation of
  reconstruction kernels.
\newblock {\em Inverse Problems}, 16(3):701, 2000.

\bibitem{Sharafutdinov_Generalized_Tensor_Fields}
Vladimir Sharafutdinov.
\newblock A problem of integral geometry for generalized tensor fields on
  {${\mathbb{R}}^n$}.
\newblock {\em Dokl. Akad. Nauk SSSR}, 286(2):305--307, 1986.

\bibitem{Sharafutdinov_Book}
Vladimir Sharafutdinov.
\newblock {\em Integral geometry of tensor fields}.
\newblock Inverse and Ill-posed Problems Series. VSP, Utrecht, 1994.

\bibitem{Sharafutdinov2007}
Vladimir Sharafutdinov.
\newblock Slice-by-slice reconstruction algorithm for vector tomography with
  incomplete data.
\newblock {\em Inverse Problems}, 23(6):2603--2627, 2007.

\bibitem{Sharafutdinov_TRT_2008}
Vladimir Sharafutdinov.
\newblock The problem of polarization tomography: II.
\newblock {\em Inverse Problems}, 24(3), 2008.

\bibitem{Sherson-15}
Brian Sherson.
\newblock {\em Some results in single-scattering tomography},
\newblock PhD Thesis, Oregon State University, 2015.


\bibitem{Sparr1995}
Gunnar Sparr, Kent Strahlen, Kjell Lindstrom and Hans W. Persson.
\newblock Doppler tomography for vector fields.
\newblock {\em Inverse problems} 11(5), pages 1051-1061, 1995.

\bibitem{Terzioglu_2015}
Fatma Terzioglu.
\newblock Some inversion formulas for the cone transform.
\newblock {\em Inverse Problems}, 31(11):115010, Oct 2015.


\bibitem{Terzioglu_2018}
Fatma Terzioglu, Peter Kuchment and Leonid Kunyansky.
\newblock Compton camera imaging and the cone transform: a brief overview.
\newblock {\em Inverse Problems}, 34(5):054002, April 2018.

\bibitem{Tuy1983}
Heang~K. Tuy.
\newblock An inversion formula for cone-beam reconstruction.
\newblock {\em SIAM Journal on Applied Mathematics}, 43(3):546--552, 1983.

\bibitem{Walker_2019} Michael R.~Walker and Joseph A.~O'Sullivan.
\newblock The broken ray transform: additional properties and new inversion formula,
\newblock \textit{Inverse Problems}, 35(11): 115003, 2019.

\bibitem{Wongsason_2018}
Patcharee Wongsason.
\newblock Vector field reconstruction via quaternionic setting.
\newblock \textit{Mathematical Methods in the Applied Sciences}, 41(2):684--696, 2018.

% \bibitem{Wongsason_Thesis}
% Patcharee Wongsason.
% \newblock {\em Reconstruction of a vector field by transverse ray transform with
%   sources on a curve}.
% \newblock PhD thesis, Oregon State University, 2014.

\bibitem{ZSM-star-14}
Fan Zhao, John C.~Schotland and Vadim A.~Markel. \newblock Inversion of the star transform, \newblock \textit{Inverse Problems}, 30(10):105001, 2014.

\end{thebibliography}
